\documentclass[12pt,twoside]{article}
\usepackage{times}
\usepackage{amsmath,amssymb}
\usepackage{color}

\allowdisplaybreaks

\def\p{\partial}
\def\ve{\varepsilon}
\def\f{\frac}
\def\na{\nabla}
\def\la{\lambda}

\def\al{\alpha}
\def\a{\alpha}

\def\t{\tilde}

\def\vp{\varphi}

\def\th{\theta}
\def\g{\gamma}
\def\G{\Gamma}
\def\si{\sigma}

\def\p{\partial}
\def\pa{\partial}
\def\ve{\varepsilon}
\def\f{\frac}

\def\na{\nabla}
\def\la{\lambda}
\def\al{\alpha}

\def\t{\tilde}

\def\o{\omega}

\def\vp{\varphi}
\def\th{\theta}

\def\g{\gamma}
\def\G{\Gamma}

\def\si{\sigma}

\def\ds{\displaystyle}
 \allowdisplaybreaks

\begin{document}
 \footskip=0pt
 \footnotesep=2pt
\let\oldsection\section
\renewcommand\section{\setcounter{equation}{0}\oldsection}
\renewcommand\thesection{\arabic{section}}
\renewcommand\theequation{\thesection.\arabic{equation}}
\newtheorem{claim}{\noindent Claim}[section]
\newtheorem{theorem}{\noindent Theorem}[section]
\newtheorem{lemma}{\noindent Lemma}[section]
\newtheorem{proposition}{\noindent Proposition}[section]
\newtheorem{definition}{\noindent Definition}[section]
\newtheorem{remark}{\noindent Remark}[section]
\newtheorem{corollary}{\noindent Corollary}[section]
\newtheorem{example}{\noindent Example}[section]

\title{The global smooth symmetric solution to 2-D
full compressible Euler system of Chaplygin gases}

\author{Ding,
Bingbing$^{1,*}$; \quad Witt, Ingo$^{2,**}$;\quad Yin,
Huicheng$^{1,}$\footnote{Ding Bingbing and Yin Huicheng are
supported by the NSFC (No.~10931007, No.~11025105), and by the
Priority Academic Program Development of Jiangsu Higher Education
Institutions.}\vspace{0.5cm}\\
\small 1. Department of Mathematics and
IMS, Nanjing University, Nanjing 210093, China.\\
\vspace{0.5cm}
\small 2.
Mathematical Institute, University of G\"{o}ttingen,
Bunsenstr.~3-5, D-37073 G\"{o}ttingen, Germany. }

\date{}
\maketitle
% \vskip 0.2in

\centerline {\bf Abstract} \vskip 0.3 true cm

For one dimensional or multidimensional compressible Euler system of polytropic gases,
it is well known that the smooth solution will generally develop singularities in finite time.
However, for three dimensional Chaplygin gases, due to the crucial role of  ``null condition'' in
the potential equation which is derived by the irrotational and isentropic flow,
P.Godin in [9] has proved the global existence of a smooth 3-D spherically symmetric flow with variable
entropy when the initial data are of small smooth perturbations with compact
supports to a constant state. It is noted that there are some essential
differences on the global solution or blowup problems between 2-D and 3-D hyperbolic systems.
In this paper, we will focus on the global symmetric solution problem of 2-D
full compressible Euler system of Chaplygin gases.
Through carrying out involved analysis and  finding an appropriate  weight
we can derive some uniform weighted energy estimates on the small
symmetric solution to 2-D compressible Euler system of Chaplygin gases and further establish the
global existence of smooth solution by continuous induction method.

\vskip 0.3 true cm

\vskip 0.3 true cm


\vskip 0.4 true cm
%\head
\centerline{\bf $\S 1$. Introduction  and main results}
%\endhead
\vskip 0.3 true cm

In this paper, we are concerned with the global existence of a smooth symmetric solution
to 2-D full compressible Euler
system of Chaplygin gases. The 2-D full Euler system is
\begin{equation}
\left\{
\begin{aligned}
&\p_t\rho+div (\rho u)=0\qquad \qquad \qquad \qquad \qquad \qquad \qquad \qquad\text{(Conservation of mass)},\\
&\p_t(\rho u)+div (\rho u \otimes u) + \nabla P=0\qquad \qquad \qquad \qquad\qquad \text{(Conservation of momentum)},\\
&\p_t\bigl(\rho ( e + \ds\frac{|u|^2}{2})\bigr)+div  \bigl( (\rho ( e + \ds\frac{|u|^2}{2}) + P
\bigr) u \bigr)=0\qquad\qquad \text{(Conservation of energy)},\\
&P=P(\rho, S),\qquad e=e(\rho, S)\qquad \qquad \qquad \qquad \qquad\qquad \text{(Equations of
state)},
\end{aligned}
\right.\tag{1.1}
\end{equation}
where $t\ge 0$, $x=(x_1, x_2)\in\Bbb R^2$, $\na=(\p_1, \p_2)$, and
$u=(u_1,u_2), \rho, P, e, S$ stand for the
velocity, density, pressure, internal energy, specific entropy
respectively. Moreover, the pressure function $P=P(\rho,S)$ and
the internal energy function $e=e(\rho,S)$ are smooth in their
arguments. In particular, $\p_{\rho}P(\rho,S)>0$ and
$\p_{S}e(\rho,S)>0$ for $\rho>0$. When $P(\rho, S)=A\rho^{\g}e^{\f{S}{c_v}}$
and $e(\rho,S)=\ds\f{A}{\g-1}\rho^{\g-1}e^{\f{S}{c_v}}$ hold for some positive
constants $A, c_v$ and $\g$ ($1<\g<3$), such flows
are then called the polytropic gases.

For the Chaplygin gases, the equation of pressure  state (one can see [7] and so on)
is given by
$$P=P_0-\ds\f{A(S)}{\rho},\eqno{(1.2)}$$
where $P_0>0$ is a positive constant, $A(S)$ is a positive smooth function of $S$,
and $P>0$ for $\rho>0$.

If $(\rho, u, S)\in C^1$ is a solution of (1.1) with $\rho>0$, then (1.1) admits
the following
equivalent form
\begin{equation}
\left\{
\begin{aligned}
&\p_t\rho+div(\rho u)=0,\\
&\p_tu+u\cdot\na u+\ds\f{\na P}{\rho}=0,\\
&\p_tS+u\cdot\na S=0.
\end{aligned}
\right.\tag{1.3}
\end{equation}
We pose the symmetric initial data of (1.3) as follows:
\begin{equation}
\left\{
\begin{aligned}
\rho(0,x)&=\bar\rho+\varepsilon\rho_0(r), \\
u(0,x)&=\varepsilon U_0(r)\ds\f{x}{r},\\
S(0,x)&=\bar S+\varepsilon S_0(r),
\end{aligned}
\right.\tag{1.4}
\end{equation}
where  $\bar\rho>0$ and $\bar S\in\Bbb R$ are constants, $\ve>0$
is a small parameter, $r=\sqrt{x_1^2+x_2^2}$ and $(\rho_0(r), U_0(r),$ $S_0(r))
\in C_0^{\infty}(B(0, M))$
(here and below $B(0, M)$ stands for a ball centered at the origin with a radius $M>0$).
Moreover, $\rho(0,x)>0$, $P_0-\ds\f{A(\bar S)}{\bar\rho}>0$ and $A'(\bar S)\not=0$ hold.

Our main result in this paper is

{\bf Theorem 1.1.} {\it Under the assumptions above, if $\ve>0$ is small enough, then (1.2)-(1.4) has a global
$C^{\infty}([0, \infty)\times\Bbb R^2)$ solution $(\rho(t,x), u(t,x), S(t,x))$ which admits such a symmetric
structure: $(\rho(t,x), u(t,x),$ $S(t,x))=(\rho(t, r), U(t,r)\ds\f{x}{r}, S(t,r))$.}

{\bf Remark 1.1.} {\it For the polytropic gases, it is well-known that the
bounded smooth solution of 1-D or multidimensional full compressible Euler system
will generally develop singularities in finite time whether the vacuum states appear or not,
one can see  [1-2], [6-8], [14], [17]
[19-23] and the references therein. However, for the Chaplygin gases, by the results of 1-D case in [17],
3-D symmetric case in [9] and our Theorem 1.1 for 2-D symmetric case, the small perturbed symmetric
flows will exist globally.
Here we point out that the  main difference for 2-D (or 3-D)
compressible Euler systems of polytropic gases and of Chaplygin gases is: if one neglects
the influences of rotations and entropy, the resulting potential equation of polytropic gases,
which is a 2-D (or 3-D) quasilinear wave equation,
does not fulfill the ``null conditions'' put forward in [4], [5] and [12]
but the related potential equation of Chaplygin gases does. In fact, when the 2-D or 3-D
quasilinear wave equations satisfy the null conditions, it is well-known that the small data
smooth solutions will exist globally (see [4-5], [11], [12] and [16]), otherwise the smooth
solution will blow up in finite time (one can see [3], [10-11] and the references therein).}

{\bf Remark 1.2.} {\it As in Remark 2.2 of [9], Theorem 1.1 still holds when
the initial data (1.4) is replaced by $(\rho(0,x), u(0,x), S(0,x))=(\bar\rho+\ve\rho_0(r),
\bar u+\ve U_0(r)\ds\f{x}{r}, \bar S+\varepsilon S_0(r))$ where $\bar u\in\Bbb R^2$ is a constant vector
since the compressible Euler system is invariable under the translation transformation.}

Let us give some comments on the proof of Theorem 1.1. As the usual first step  to prove the global existence or blowup
of a smooth small data solution for the quasilinear hyperbolic equation, one should construct a suitable approximate solution
$(\rho_a, u_a, S_a)$ to (1.2)-(1.4). To this end, as in [8-9] and [21], a solution to the
related potential equation with suitable initial data will be taken as the approximate solution.
Through considering the difference of the real solution $(\rho, u, S)$ and approximate solution
$(\rho_a, u_a, S_a)$, and by applying a transformation of unknowns introduced in [21],
one can obtain a new symmetric hyperbolic system.  From this system, we can get a potential
equation satisfying  both null conditions in two space dimensions outside a compact support of
perturbed entropy. The so-called both null conditions in two space dimensions are defined in [4] as follows:
For the 2-D quasilinear wave equation $\square w+\ds\sum_{0\le i,j,k\le 2}g_{ij}^k\p_kw\p_{ij}^2w+
\ds\sum_{0\le i,j,k,l\le 2}g_{ij}^{kl}\p_kw\p_lw\p_{ij}^2w=0$, where $x_0=t$, $g_{ij}^k$ and
$g_{ij}^{kl}$ are certain constants, then $\ds\sum_{0\le i,j,k\le 2}g_{ij}^k\xi_k\xi_i\xi_j\equiv 0$
and $\ds\sum_{0\le i,j,k,l\le 2}g_{ij}^{kl}\xi_k\xi_l\xi_i\xi_j\equiv 0$  hold for $(\xi_0, \xi_1, \xi_2)
=(-1, cos\th, sin\th)$ with $0\le\th\le 2\pi$. Especially, the former
$\ds\sum_{0\le i,j,k\le 2}g_{ij}^k\xi_k\xi_i\xi_j\equiv 0$ and the
latter $\ds\sum_{0\le i,j,k,l\le 2}g_{ij}^{kl}\xi_k\xi_l\xi_i\xi_j\equiv 0$ are called the first and the second null condition respectively
by the terminology in [4]. Under both null conditions, S. Alinhac in [4] established the global existence of small data solution
to the 2-D quasilinear wave equation by looking for a crucial ``ghost weight"  to derive an energy estimate.
In this paper, we will focus on the global existence of solution to (1.2)-(1.4). To this end,
we require to derive the uniform weighted energy estimate of solution $(\rho, u, S)$ so that the uniform bound
of $(\rho, u, S)$  for $(t,x)\in [0, \infty)\times\Bbb R^2$ can be obtained. Motivated by the methods in [4],
we will look for a suitable ``ghost weight'' to deal with the related compressible  Euler system
(1.3) so that both null conditions and the variable entropy can be simultaneously considered.
Here we point out that the ghost weight introduced in [4] will not be applied for our case directly
since (1.3) can not be changed into a scalar quasilinear wave equation due to the influence of variable entropy,
and thus the ghost weight in [4] should be suitably adjusted for our uses (one can see more detailed explanations
in Remark 4.1 of $\S 4$ below).
On the other hand, although some procedures in this paper are somewhat analogous to those in [8-9] for considering the 3-D
symmetric Euler system, our analysis
is more involved since the decay rate of solution to 2-D free wave equation is lower
than that in 3-D case as well as the treatments on the both null conditions in 2-D quasilinear wave equation
are more complicated than the treatments on one null condition in 3-D quasilinear wave equation (one can compare the reference
[4] with [5] and [12]).

The rest of the paper is organized as follows. In $\S 2$, we will
construct an approximate solution of (1.2)-(1.4) and give some useful preliminary
knowledge. In $\S 3$ and $\S 4$, we will establish the uniform higher order
energy estimates on the solution $(\rho, u, S)$ near the light cone and in the whole time-space respectively.
Based on these uniform estimates, Theorem 1.1 will be proved by continuous induction method.

In what follows, we will use the following convention as in [8-9] and [21]:

$x_0=t$ $(\geq 0)$ denotes the time variable, $\partial=(\partial_0,\partial_1,\partial_2)=(\partial_t,\nabla)$;

$
\Omega=x^\perp\cdot\nabla$ with $x^\perp=(-x_2,x_1)$, $
X=t\partial_t+r\partial_r
$, $
L_j=t\partial_j+x_j\partial_t$ for $j=1,2$;

$\Lambda=(\Lambda_1, \cdots, \Lambda_7)=(\partial, \Omega, X, L_1, L_2),
$ $\Gamma=(\Gamma_1,\cdots,\Gamma_4)=(\partial, X)$;

$Z=(Z_1,Z_2)=(\pa_1+\omega_1\pa_t,\pa_2+\omega_2\pa_t)=\ds\frac{t-r}{t}(\pa_1,\pa_2)
+\frac{(x_1,x_2)X+(-x_2,x_1)\Omega}{rt}$ with $\o_i=\ds\f{x_i}{r}$ ($i=1,2$) and $\o=(\o_1, \o_2)$;

For $\Bbb{R}^l$-valued smooth functions $f(t,x)$ and $\tilde f(t,x)$, we put
$$
|f(t)|=\sup_{x\in\Bbb{R}^2}|f(t,x)|,\quad |f(t)|_\pm=\sup_{D_\pm}|f(t,x)|,
$$
$$
\langle f,\tilde f\rangle(t)=\int_{\Bbb{R}^2}f(t,x)\cdot\tilde f(t,x)\,{d x},
\quad \langle f,\tilde f\rangle_\pm(t)=\int_{D_\pm(t)}f(t,x)\cdot\tilde f(t,x)\,{d x},
$$
and
$$
\|f(t)\|=\langle f,f\rangle^{1/2}(t),\quad \|f(t)\|_\pm=\langle f,f\rangle_\pm^{1/2}(t),
$$
where $D_-(t)=\{x\in\Bbb{R}^2:
r\leq\ds\frac{t}{2}+M+1\}$ and $D_+(t)=\{x\in\Bbb{R}^2:
r\geq\ds\frac{t}{2}+M+1\}$;

If $g$ is only a function of $x$, then $\|g\|$  represents the usual $L^2$ norm;

Define $\langle\xi\rangle=1+|\xi|$ for $\xi\in\Bbb R$ and $\sigma(x)=\langle x\rangle$,
$\sigma_\pm(t,x)=\langle t\pm |x|\rangle$.

\vskip 0.4 true cm \centerline{\bf $\S 2$. The construction of approximate solution to
(1.2)-(1.4) and some preliminaries} \vskip 0.4 true cm

Without loss of generality, we will assume the sound speed
$\bar c=\big(\pa_{\rho} P(\bar\rho,\bar S)\big)^{1/2}\equiv 1$
in the whole paper. As in [8-9], denote $(\rho_a, u_a, \bar S)$ by a solution of
(1.3)-(1.4) satisfying the initial data condition $u_a(0,x)=u_0(0,x)$ and $P_a(0,x)=P(0,x)$, which means
the initial density $\rho_a(0,x)=\rho_a(0,r)=\ds\frac{A(\bar S)\rho(0,x)}{A(S(0,x))}$. At this time, since the initial
data of $(\rho_a, u_a, \bar S)$ are symmetric and isentropic, one can introduce a potential $\phi_a(t,r)$
such that $u_a=\na\phi_a$ and $\phi_a$ satisfies the following potential equation
\begin{equation}
\left\{
\begin{aligned}
&\pa_t^2\phi_a+2\ds\sum_{j=1}^2\pa_j\phi_a\pa_t\pa_j\phi_a+\ds\sum_{j,k=1}^2\pa_j\phi_a\pa_k\phi_a\pa_{jk}^2\phi_a
-(1+2\pa_t\phi_a+|\nabla \phi_a|^2)\triangle \phi_a=0,\\
&\phi_a(0,r)=-\displaystyle\ve\int_{r}^MU_0(s)ds,\\
&\pa_t\phi_a(0,r)=-\displaystyle\frac{1}{2}\ve^2U_0^2(r)-h\big(\rho_a(0,r)\big),
\end{aligned}
\right.\tag{2.1}
\end{equation}
where $h(\rho)=\ds\f12-\ds\f{A(\bar S)}{2\rho^2}$ is the enthalpy. Meanwhile, the density
$\rho_a$ is determined by the Bernoulli's law $h(\rho_a)=-\p_t\phi_a-\f12|\na\phi_a|^2$.

It is easy to verify that (2.1) satisfies both null conditions in two space dimensions
(i.e., the first null condition and second null condition, which have been illustrated in $\S 1$) posed in [4]. Then
(2.1) has a global smooth solution $\phi_a$ in terms of [4] and we have

{\bf Lemma 2.1.}
{\it If $\ve>0$ is small enough, then

$(1)\quad(|\Lambda^{\al}(\rho_a-\bar\rho)|+|\Lambda^{\al}u_a|)(t,x)
\leq C_{\al}\ve\si_-(t,x)^{-1}\si_+(t,x)^{-1/2}$;

$(2)\quad|\pa\Lambda^{\al}(\rho_a-\bar\rho)(t,x)|+|\pa\Lambda^{\al}u_a(t,x)|
\leq C_{\al}\ve\langle t\rangle^{-5/2}$  if $|x|\leq lt+M$ with $0\leq l<1$;

$(3)\quad|X^ku_a(t,x)|\leq C_{k\delta}\ve\langle t\rangle^{-5/2+\delta}$ if $0\leq\delta<1$
and $|x|\leq C\langle t\rangle^{\delta}$,

where $C_{\al}$, $C_{k\delta}$ and $C$ are some generic positive constants independent of $\ve$ and $(t,x)$.}

{\bf Proof.} Since (2.1) satisfies both null conditions in 2-D spaces, then one
has by the result of [4] or [13]
$$
|\Lambda^{\al}\pa^\beta \phi_a|\leq C_{\al}\ve\si_-(t,x)^{-|\beta|}\si_+(t,x)^{-1/2}.\eqno{(2.2)}
$$
From this, (1) and (2) can be obtained directly.

In addition, by $\ds X^ku_a(t,x)
=(X^k\p_r\phi_a)(t,r)\frac{x}{r}
=\biggl(\int_0^r\pa_\lambda (X^k\p_r\phi_a)(t,\lambda)\,{d \lambda}\biggr)\frac{x}{r}$ and (2.2), then (3) holds
obviously.
\qquad\qquad\qquad\qquad  \qquad\qquad  \qquad\qquad  \qquad\qquad  \qquad\qquad     $\square$\bigskip

As in [8] and [21], we set $\theta(t,x)=1-\ds\frac{A(S(t,x))\bar\rho}{A(\bar S)\rho(t,x)}$,
$w(t,x)=u(t,x)$, $z(t,x)=\ds\frac{A(\bar S)}{A(S(t,x))}-1$, then it follows from (1.3)-(1.4) that
\begin{equation}
\left\{
\begin{aligned}
&\partial_t\theta+w\cdot\nabla\theta+(1-\theta)\nabla\cdot w=0,\\
&\partial_tw+w\cdot\nabla w+(1-\theta)(1+z)\nabla \theta=0,\\
&\partial_tz+w\cdot\nabla z=0,\\
&\theta(0,x)=1-\ds\frac{A(\bar S+\ve S_0(r))\bar\rho}{A(\bar S)(\bar\rho+\ve\rho_0(r))},\\
&w(0,x)=\varepsilon U_0(r)\ds\f{x}{r},\\
&z(0,x)=\ds\frac{A(\bar S)}{A(\bar S+\ve S_0(r))}-1.
\end{aligned}
\right.\tag{2.3}
\end{equation}

Corresponding to the approximate solution $(\rho_a,u_a,\bar S)$, define $\theta_a(t,x)=1-\ds\frac{\bar\rho}{\rho_a(t,x)}$
and $w_a(t,x)=u_a(t,x)$. Then a direct computation yields by Lemma 2.1 that

{\bf Lemma 2.2.}
{\it For small $\ve>0$, we have

$(1)\quad(|\Lambda^{\al}\theta_a|+|\Lambda^{\al}w_a|)(t,x)\leq C_{\al}\ve\si_-(t,x)^{-1}\si_+(t,x)^{-1/2}$;

$(2)\quad|\pa\Lambda^{\al}\theta_a(t,x)|+|\pa\Lambda^{\al}w_a(t,x)|\leq C_{\al}\ve\langle t\rangle^{-5/2}$
if  $|x|\le lt+M$ with $0\le l< 1$;

$(3)\quad|X^kw_a(t,x)|\leq C_{k\delta}\ve\langle t\rangle^{-5/2+\delta}$ if $0\leq\delta<1$ and $|x|\leq C\langle t\rangle^{\delta}$.\hfil\break}

Let $\dot\theta=\theta-\theta_a$, $\dot w=w-w_a$ and $\dot z=z-z(0,x)$, then we have by (2.3)
\begin{equation}
\left\{
\begin{aligned}
&\partial_t\dot\theta+\nabla\cdot \dot w=-(w\cdot\nabla\theta-w_a\cdot\nabla\theta_a)+(\theta\nabla\cdot w
-\theta_a\nabla\cdot w_a),\\
&\partial_t\dot w+\nabla\dot \theta=-(w\cdot\nabla w-w_a\cdot\nabla w_a)+(\theta\nabla\theta-\theta_a\nabla\theta_a)-(1-\theta)z\nabla\theta,\\
&\partial_t\dot z=-w\cdot\nabla(z(0,x)+\dot z), \\
&\dot \theta(0,x)=0,\quad\dot w(0,x)=0,\quad \dot z(0,x)=0.
\end{aligned}
\right.\tag{2.4}
\end{equation}

To establish the global existence of solution to (2.4) in subsequent sections, we require to
give some preliminary analysis
on the related energies.
As usual, we define the energy for $n\in\Bbb N\cup\{0\}$
$$
E_n(t)=\sum_{|\al|\leq n}(\|\Gamma^{\al}\dot \theta(t)\|^2
+\|\Gamma^{\al} \dot w(t)\|^2+\|\Gamma^{\al} \dot z(t)\|^2).
$$

We also set $Q_n(t)=\ds\sum_{|\al|\leq n-1}\bigg(\|\si_-(t)\nabla\Gamma^{\al}\dot \theta(t)\|
+\|\si_-(t)\pa_t\Gamma^{\al}\dot w(t)\|
+\|\si_-(t)\nabla\cdot\Gamma^{\al}\dot w(t)\|\bigg)$ for $n\geq 1$, and define
$\widetilde Q_n(t)=Q_n(t)+E_{n-1}^{1/2}(t)$
for $n\geq 1$ and $\widehat Q_n(t)=Q_n(t)+E_{n-2}^{1/2}(t)$ for $n\geq 2$ as in [8-9]
and [21].

In addition, for our requirements to treat the 2-D full Euler system (1.3)-(1.4), it is necessary
to introduce some kinds of ``interior energies'' as follows:

Choose a smooth function
\begin{equation*}
\hat\chi(s) = \left\{
\begin{aligned}
& 1,  \quad \text{if $s\leq1/2$},\\
& 0, \quad \text{if $s\geq3/4$},
\end{aligned}
\right.
\end{equation*}
and set $\chi(t,x)=\hat\chi(\ds\frac{|x|}{t+2M+2})$, then we define for $n\geq 1$
\begin{align*}
&Q_n^-(t)=\ds\sum_{|\al|\leq n-1}(\|\si_-(t)\nabla\Gamma^{\al}(\chi\dot\theta)(t)\|
+\|\si_-(t)\pa_t\Gamma^{\al}(\chi \dot w)(t)\|
+\|\si_-(t)\nabla\cdot\Gamma^{\al}(\chi \dot w)(t)\|),\\
&\widetilde Q_n^-(t)=Q_n^-(t)+E_{n-1}^{1/2}(t),\quad\text{for $n\geq 1$},\\
&\widehat Q_n^-(t)=Q_n^-(t)+E_{n-2}^{1/2}(t),\quad \text{for $n\geq 2$}.
\end{align*}

The following relations and properties on the above defined energies will be repeatedly utilized
later on.

{\bf Lemma 2.3.} {\it
\begin{align*}
(1)\quad&\|\si_-(t)\nabla v\|\leq C(\|\si_-(t)\nabla\cdot v\|+\|v\|)\quad\text{if $v\in C_0^\infty(\Bbb R^2,\Bbb R^2)$ and $\nabla\cdot v^\bot=0$,}\\
(2)\quad&|\si^{1/2}\Gamma^{\al}\dot w(t)|+|\si^{1/2}\Gamma^{\al}\dot\theta(t)|
+|\si^{1/2}\Gamma^{\al}\dot z(t)|\leq C_{\al}E_{|\al|+2}^{1/2}(t),\\
(3)\quad&|\si^{1/2}\si_-(t)\nabla\Gamma^{\al}\dot\theta(t)|\leq C_{\al}Q_{|\al|+3}(t),\\
(4)\quad&|\si^{1/2}\si_-(t)\nabla\Gamma^{\al}(\chi\dot\theta)(t)|\leq C_{\al}Q^-_{|\al|+3}(t),\\
(5)\quad&|\si^{1/2}\si_-(t)\nabla\Gamma^{\al}\dot w(t)|\leq C_{\al}\widetilde Q_{|\al|+3}(t),\\
(6)\quad&|\si^{1/2}\si_-(t)\nabla\Gamma^{\al}(\chi \dot w)(t)|\leq C_{\al}\widetilde Q^-_{|\al|+3}(t),\\
(7)\quad&|\si^{1/2}\si_-^{1/2}(t)\Gamma^{\al}\dot\theta(t)|\leq C_{\al}\widehat Q_{|\al|+2}(t),\\
(8)\quad&|\si^{1/2}\si_-^{1/2}(t)\Gamma^{\al}(\chi\dot\theta)(t)|\leq C_{\al}\widehat Q^-_{|\al|+2}(t),\\
(9)\quad&|\si^{1/2}\si_-^{1/2}(t)\Gamma^{\al}\dot w(t)|\leq C_{\al}\widetilde Q_{|\al|+2}(t),\\
(10)\quad&|\si^{1/2}\si_-^{1/2}(t)\Gamma^{\al}(\chi \dot w)(t)|\leq C_{\al}\widetilde Q_{|\al|+2}^-(t).
\end{align*}
}

{\bf Proof.} (1) comes from the formula (6.7) in [21] directly.

In addition, according to Lemma 1 of [21], one has for any smooth function $f(x)$
\begin{equation}
\left\{
\begin{aligned}
&\si(x)^{1/2}|f(x)|\leq C\ds\sum_{j=0}^1\sum_{|\a|=0}^{2-j}\|\nabla^\a\Omega^jf\|,\\
&\si(x)^{1/2}\si_-(t,x)|f(x)|\leq C\ds\sum_{j=0}^1\sum_{|\a|=0}^{2-j}\|\si_-(t)\nabla^\a\Omega^jf\|,\\
&\si(x)^{1/2}\si_-(t,x)^{1/2}|f(x)|\leq C\ds\sum_{j=0}^1\big(\|\Omega^jf\|+\sum_{|\a|=1}^{2-j}\|\si_-(t)\nabla^\a\Omega^jf\|\big).
\end{aligned}
\right.\tag{2.5}
\end{equation}

This, together with the facts of $\nabla\cdot(\Gamma^a(\chi \dot w))^\bot=0$
and $[\Omega,\Gamma^{\al}]=\ds\sum_{|\beta|\leq|\al|}C_{\al\beta}\Gamma^{\beta}$, yields
$(2)-(10)$ directly.\qquad \qquad \qquad \qquad \qquad \qquad \qquad $\square$
\bigskip

Applying $\pa^\a(X+1)^k$ to (2.4) and setting $\Gamma^{\mu}=\pa^\a X^k$, then we have
\begin{equation}
\left\{
\begin{aligned}
\pa_t\Gamma^{\mu}\dot\theta+\nabla\cdot\Gamma^{\mu}\dot w&=h_0^{\mu}\equiv \ds\sum_{1\leq j\leq6}f_j^{\mu},\\
\pa_t\Gamma^{\mu}\dot w+\nabla\Gamma^{\mu}\dot\theta&=h^{\mu}\equiv \ds\sum_{7\leq j\leq 13}f_j^{\mu},
\end{aligned}
\right.\tag{2.6}
\end{equation}
\break
$f_1^{\mu}=-\ds\sum_{{\nu}\leq \mu}\binom{\mu}{\nu}\Gamma^{\nu}w_a\cdot\nabla\Gamma^{\mu-\nu}\dot\theta,
\quad f_2^{\mu}=-\ds\sum_{\nu\leq \mu}\binom{\mu}{\nu}\Gamma^{\nu}\dot w\cdot\nabla\Gamma^{\mu-\nu}\theta_a,$\hfil\break
$f_3^{\mu}=-\ds\sum_{\nu\leq \mu}\binom{\mu}{\nu}\Gamma^{\nu}\dot w\cdot\nabla\Gamma^{\mu-\nu}\dot\theta,\quad
f_4^{\mu}=\ds\sum_{\nu\leq \mu}\binom{\mu}{\nu}\Gamma^{\nu}\theta_a\nabla\cdot\Gamma^{\mu-\nu}\dot w,$\hfil\break
$f_5^{\mu}=\sum_{\nu\leq \mu}\binom{\mu}{\nu}\Gamma^{\nu}\dot\theta\nabla\cdot\Gamma^{\mu-\nu}w_a,\quad
f_6^{\mu}=\sum_{\nu\leq \mu}\binom{\mu}{\nu}\Gamma^{\nu}\dot\theta\nabla\cdot\Gamma^{\mu-\nu}\dot w,$\hfil\break
$f_7^{\mu}=-\ds\sum_{\nu\leq\mu}\binom{\mu}{\nu}\Gamma^{\nu}w_a\cdot\nabla\Gamma^{\mu-\nu}\dot w,\quad f_8^{\mu}=-\ds\sum_{\nu\leq\mu}\binom{\mu}{\nu}\Gamma^{\nu}\dot w\cdot\nabla\Gamma^{\mu-\nu}w_a,
$\hfil\break
$f_9^{\mu}=-\ds\sum_{\nu\leq\mu}\binom{\mu}{\nu}\Gamma^{\nu}\dot w\cdot\nabla\Gamma^{\mu-\nu}\dot w,
\quad f_{10}^{\mu}=\ds\sum_{\nu\leq \mu}\binom{\mu}{\nu}\Gamma^{\nu}\theta_a\nabla\Gamma^{\mu-\nu}\dot\theta,$\hfil\break
$f_{11}^{\mu}=\sum_{\nu\leq \mu}\binom{\mu}{\nu}\Gamma^{\nu}\dot\theta\nabla\Gamma^{\mu-\nu}\theta_a,\quad
f_{12}^{\mu}=\sum_{\nu\leq\mu}\binom{\mu}{\nu}\Gamma^{\nu}\dot\theta\nabla\Gamma^{\mu-\nu}\dot\theta,$\hfil\break
$f_{13}^{\mu}=-\pa^\a(X+1)^k\big((1-\theta)z\nabla\theta\big)$,
\quad here we point out that the concrete expressions of
$f_i^{\mu}$ ($1\le i\le 12$) are important in order to derive the basic energy estimate in Lemma 3.5 below
(one can see the details in dealing with (3.21)) since
we require to utilize the first null condition by observing the main parts of $f_1^{\mu}$ and $f_4^{\mu}$ (or $f_7^{\mu}$ and $f_{10}^{\mu}$), and the second null condition by observing the left parts of $f_1^{\mu}$ and $f_4^{\mu}$
(or $f_7^{\mu}$ and $f_{10}^{\mu}$) except the terms in the first null condition together with the main parts of $f_i^{\mu}$  for $i= 2, 3, 5, 6$ (or $f_i^{\mu}$  for $i=8, 9, 11, 12$).

It follows from an analogous computation in [9] that
\begin{align*}
&(t-r)\nabla\Gamma^{\mu}\dot\theta=-\frac{t}{t+r}X\Gamma^{\mu}\dot w+\frac{t}{t+r}(\Omega\Gamma^{\mu}\dot w)^\bot-(X\Gamma^{\mu}\dot\theta)\frac{x}{t+r}\\
&\qquad \qquad \qquad -\frac{x^\bot}{t+r}\Omega\Gamma^{\mu}\dot\theta
+\frac{t}{t+r}h_0^{\mu}x+\frac{t^2}{t+r}h^{\mu},\tag{2.7}\\
&(t-r)\nabla\cdot\Gamma^{\mu}\dot w=-\frac{t}{t+r}X\Gamma^{\mu}\dot\theta-\frac{X\Gamma^{\mu}\dot w}{t+r}\cdot x-\frac{(\Omega\Gamma^{\mu}\dot w)^\bot}{t+r}\cdot x+\frac{t^2}{t+r}h_0^{\mu}
+\frac{t}{t+r}h^{\mu}\cdot x,\tag{2.8}\\
&(t-r)\pa_t\Gamma^{\mu}\dot w=\frac{t}{t+r}X\Gamma^{\mu}\dot w-\frac{t}{t+r}(\Omega\Gamma^{\mu}\dot w)^\bot+(X\Gamma^{\mu}\dot\theta)\frac{x}{t+r}+\frac{x^\bot}{t+r}\Omega\Gamma^{\mu}\dot\theta\\
&\qquad \qquad \qquad -\frac{t}{t+r}h_0^{\mu}x-\frac{r^2}{t+r}h^{\mu}.\tag{2.9}
\end{align*}

Based on (2.7)-(2.9), we can establish a useful energy inequality inside the light cone as follows:

{\bf Lemma 2.4.} {\it
For fixed $n\in\Bbb N$ with $n\geq 4$. There exist positive constants $\eta$ and $C$ such that
if $(\dot\theta, \dot w, \dot z)$ is a smooth solution of (2.4)
for $(t,x)\in [0, T]\times\Bbb R^2$ with $E_{[\frac{n}{2}]+2}^{1/2}(t)\leq\eta$, then for $0\leq t\leq T$ and small $\ve>0$
$$\ds Q_n^-(t)\leq C\biggl(E_n^{1/2}(t)+\frac{\ve^2}{\langle t\rangle^{3/2}}\biggr).\eqno{(2.10)}$$
}
\bigskip

{\bf Remark 2.1.} {\it For the free wave equation with compactly supported initial data,
it is well-known that the decay rate of smooth solution on the time $t$
in 2-D case is slower than that in 3-D case, therefore the author in [9] can obtain
an energy estimate of $Q_n(t)$
similar to (2.10) in the whole space $\Bbb R^3$ (one can see Proposition 4.1 of [9]) but
at present we only get (2.10) for
$Q_n^-(t)$ which is a kind of interior energy.}

{\bf Proof.} By (2.6), we have for $|\mu|\le n-1$
\begin{equation}
\left\{
\begin{aligned}
&\pa_t\Gamma^{\mu}(\chi\dot\theta)+\nabla\cdot\Gamma^{\mu}(\chi \dot w)=\tilde h_0^{\mu}\equiv
\ds\sum_{0\leq j\leq 6}\tilde f_j^{\mu},\\
&\pa_t\Gamma^{\mu}(\chi\dot w)+\nabla\Gamma^{\mu}(\chi\dot\theta)=\tilde h^{\mu}\equiv \ds\sum_{7\leq j\leq 14}\tilde f_j^{\mu},
\end{aligned}
\right.\tag{2.11}
\end{equation}
where\hfil\break
$\tilde f_0^{\mu}=\pa^\a(X+1)^k\bigg\{\ds\frac{\chi'}{t+2M+2}\biggl(-\frac{r}{t+2M+2}\dot\theta
+\dot w\cdot\omega+w\cdot\omega\dot\theta-\theta\omega\cdot \dot w\biggr)\bigg\},$\hfil\break
$\tilde f_3^\mu=-\ds\sum_{\nu\leq\mu}\binom{\mu}{\nu}\Gamma^\nu\dot w\cdot\nabla\Gamma^{\mu-\nu}(\chi\dot\theta),\quad\tilde f_6^\mu=\ds\sum_{\nu\leq\mu}\binom{\mu}{\nu}\Gamma^\nu\dot\theta\nabla\cdot\Gamma^{\mu-\nu}(\chi\dot w),
$\hfil\break
$\tilde f_9^\mu=-\ds\sum_{\nu\leq\mu}\binom{\mu}{\nu}\Gamma^\nu\dot w\cdot\nabla\Gamma^{\mu-\nu}(\chi \dot w),\quad
\tilde f_{12}^\mu=\ds\sum_{\nu\leq\mu}\binom{\mu}{\nu}\Gamma^\nu\dot\theta\nabla\Gamma^{\mu-\nu}(\chi\dot\theta),
$\hfil\break
$\ds\tilde f_{13}^\mu=-\pa^\a(X+1)^k\big(\chi(1-\theta)z\nabla\theta\big),$\hfil\break
$\ds\tilde f_{14}^\mu=\pa^\a(X+1)^k\bigg\{\frac{\chi'}{t+2M+2}\bigg(-\frac{r}{t+2M+2}
\dot w+\dot\theta\omega+w\cdot\omega \dot w-\theta\dot\theta\omega\bigg)\bigg\},$\hfil\break
and for $i=1,2,4,5,7,8,10,11$, the expressions of $\t f_i^{\mu}$ are the same as $f_i^{\mu}$ in (2.6)
when $\dot\th$ or $\dot w$
in $f_i^{\mu}$ is replaced by $\chi\dot\th$ or $\chi\dot w$ respectively.

By (2.11) together with the similar expressions of (2.7)-(2.9), an easy computation yields

$$
Q_n^-(t)\leq C\big(E_n^{1/2}(t)+\sum_{|\mu|\leq n-1}t\|\tilde h_0^\mu(t)\|+\sum_{|\mu|\leq n-1}\langle t\rangle\|\tilde h^\mu(t)\|\big).\eqno{(2.12)}
$$

We now focus on the estimates of $\|\tilde h_0^\mu(t)\|$ and $\|\tilde h^\mu(t)\|$ in the right hand side of (2.12).

According to Lemma 2.2 and by direct observations, we can obtain
$$
\|\tilde f_j^\mu\|\leq C_n\ve\langle t\rangle^{-3/2} E_n^{1/2}(t),\quad j\in\{1,2,4,5,7,8,10,11\},\eqno{(2.13)}
$$
and
$$
\|\tilde f_0^\mu(t)\|+\|\tilde f_{14}^\mu(t)\|\leq C_n\langle t\rangle^{-1}E_{n-1}^{1/2}(t).\eqno{(2.14)}
$$

In addition, applying Lemma 2.3 (2) and (4) respectively yield that

if $|\nu|\leq|\mu-\nu|$, then
$$
\|\Gamma^\nu\dot w\cdot\nabla\Gamma^{\mu-\nu}(\chi\dot\theta)(t)\|\leq C\langle t\rangle^{-1}|\Gamma^\nu \dot w(t)|\cdot\|\si_-(t)\nabla\Gamma^{\mu-\nu}(\chi\dot\theta)(t)\|\leq C_n\langle t\rangle^{-1}E_{[\frac{n-1}{2}]+2}^{1/2}(t)Q_n^-(t);
$$

if $|\nu|>|\mu-\nu|$, then
$$
\|\Gamma^\nu\dot w\cdot\nabla\Gamma^{\mu-\nu}(\chi\dot\theta)(t)\|\leq |\nabla\Gamma^{\mu-\nu}(\chi\dot\theta)(t)|\cdot\|\Gamma^{\nu}\dot w(t)\|\leq C_n\langle t\rangle^{-1}Q^-_{[\frac{n}{2}]+2}(t)E^{1/2}_{n-1}(t).
$$

Therefore, one has
$$\|\tilde f_3^\mu\|\leq C_n\langle t\rangle^{-1}\big(E_{[\frac{n-1}{2}]+2}^{1/2}(t)Q_n^-(t)+Q^-_{[\frac{n}{2}]+2}(t)E^{1/2}_{n-1}(t)\big).\eqno{(2.15)}$$

Similarly,
$$
\|\tilde f_j^\mu\|\leq C_n\langle t\rangle^{-1} \big(E_{[\frac{n-1}{2}]+2}^{1/2}(t) \widetilde Q_n^-(t)+\widetilde Q^-_{[\frac{n}{2}]+2}(t)E^{1/2}_{n-1}(t)\big),\quad j\in\{6,9,12\}.\eqno{(2.16)}
$$

Finally, we deal with $\tilde f_{13}^\mu$.

Set
\begin{align*}
&{L}_{1\nu}(t)=\|\Gamma^\nu z(0,x)\nabla\Gamma^{\mu-\nu}(\chi\theta_a)(t)\|,\\
&{L}_{2\nu}(t)=\|\Gamma^\nu z(0,x)\nabla\Gamma^{\mu-\nu}(\chi\dot\theta)(t)\|,\\
&{L}_{3\nu}(t)=\|\Gamma^\nu\dot z\nabla\Gamma^{\mu-\nu}(\chi\theta_a)(t)\|+\|\pa^\a(X+1)^k(\frac{\chi'}{t+2M+2}\dot z\theta_a\omega)(t)\|,\\
&{L}_{4\nu}(t)=\|\Gamma^\nu\dot z\nabla\Gamma^{\mu-\nu}(\chi\dot\theta)(t)\|+\|\pa^\a(X+1)^k(\ds\frac{\chi'}{t+2M+2}\dot z\dot\theta\omega)(t)\|.
\end{align*}

It is easy to obtain
\begin{align*}
&{L}_{1\nu}(t)\le C_n\ve^2\langle t\rangle^{-5/2},\\
&{L}_{2\nu}(t)\le C_n\ve\langle t\rangle^{-1}Q_n^-(t),\\
&{L}_{3\nu}(t)\leq C_n\ve\langle t\rangle^{-5/2}E_{n-1}^{1/2}(t),\\
&{L}_{4\nu}(t)\leq C_n\langle t\rangle^{-1}\bigg(E_{[\frac{n-1}{2}]+2}^{1/2}(t)Q_n^-(t)
+Q^-_{[\frac{n}{2}]+2}(t)E^{1/2}_{n-1}(t)\bigg)
+C_n\langle t\rangle^{-3/2}E_{[\frac{n-1}{2}]+2}^{1/2}(t)E_{n-1}^{1/2}(t).
\end{align*}

This, together with the expression of $\tilde f_{13}^\mu$, yields
$$\|\tilde f_{13}^\mu\|\leq C_n\big(1+E_{[\frac{n-1}{2}]+2}^{1/2}(t)\big)F_{n-1}(t)+C_nE_{n-1}^{1/2}(t)F_{[\frac{n}{2}]+1}(t),
\eqno{(2.17)}$$
where $F_{j-1}(t)$ is defined as $F_{j-1}(t)=\ve^2\langle t\rangle^{-5/2}+E_{j-1}^{1/2}(t)\big(\ve\langle t\rangle^{-5/2}+\langle t\rangle^{-1}Q^-_{[\frac{j}{2}]+2}(t)\big)+\langle t\rangle^{-1}Q_j^-(t)\big(\ve+E_{[\frac{j-1}{2}]+2}^{1/2}(t)\big)+\langle t\rangle^{-3/2}E_{[\frac{j-1}{2}]+2}^{1/2}(t)E_{j-1}^{1/2}(t)$ for $j\ge 1$.

Due to $E^{\f12}_{[\frac{n}{2}]+2}(t)\leq\eta$, (2.12) together with (2.13)-(2.17) derives
$$
Q_{[\frac{n}{2}]+2}^-(t)\leq C_n\eta\big(1+Q_{[\frac{n}{2}]+2}^-(t)\big)+C_n\ve^2\langle t\rangle^{-3/2},
$$
which means for small $\eta$
$$
Q_{[\frac{n}{2}]+2}^-(t)\leq C_n\big(\eta+\frac{\ve^2}{\langle t\rangle^{3/2}}\big).\eqno{(2.18)}
$$
And hence, $F_{n-1}(t)\leq C_n\langle t\rangle^{-1}\bigg(Q_n^-(t)(\ve+\eta)+\eta E_{n-1}^{1/2}(t)+\ds\frac{\ve^2}{\langle t\rangle^{3/2}}\bigg)$ and $F_{[\frac{n}{2}]+1}(t)\leq C_n\langle t\rangle^{-1}\bigg(\eta(\ve+\eta)+\ds\frac{\ve^2}{\langle t\rangle^{3/2}}\bigg)$. Substituting this into (2.17) and further combining with (2.12)-(2.16) yield
$$
Q_n^-(t)\leq C_n\bigg(E_n^{1/2}(t)+\frac{\ve^2}{\langle t\rangle^{3/2}}\bigg).\qquad\qquad\qquad\qquad\text{$\square$}
$$

Next, we derive the decay estimate of solution $(\dot\th,\dot w)$ inside the light cone, which is analogous to Proposition 4.2 in [9]. This is relatively easier
to be obtained than the one near the cone in subsequent $\S 3$.

{\bf Lemma 2.5.} {\it
For fixed $\la\in\Bbb N$ with $\la\geq 4$, if $(\dot\theta, \dot w, \dot z)$ is a smooth solution of (2.4) for $(t,x)\in [0, T]\times\Bbb R^2$ with $E_\la(t)\leq\ve^2$, and  $\widetilde Q_{|\mu|+2}(t)\leq\eta\langle t\rangle^{1/2}$ holds for $|\mu|\geq\la-1$ and small constant $\eta>0$, then\hfil\break
$(1)\quad|\nabla\Gamma^\mu(\chi\dot\theta)(t)|+|\nabla\cdot\Gamma^\mu(\chi \dot w)(t)|\leq C\chi_{|\mu|}(t)$ for $|\mu|\leq 2\la-4$, where and below $\chi_k(t)\equiv \ds\frac{\widetilde Q_{k+3}^-(t)}{\langle t\rangle^{3/2}}+\frac{\ve^2}{\langle t\rangle^{5/2}}$.\hfil\break
$(2)\quad|\pa_t\Gamma^\mu(\chi\dot\theta)(t)|+|\pa_t\Gamma^\mu(\chi\dot w)(t)|\leq C\chi_{|\mu|}(t)$ for $|\mu|\leq 2\la-4$.
\hfil\break
$(3)\quad|\Gamma^\mu(\chi\dot w)(t)|\leq C\chi_{|\mu|-1}(t)$ for $|\mu|\leq 2\la-3$ and $\mu_2+\mu_3\neq0$.\hfil\break
$(4)\quad|\Gamma^\mu\dot w(t,x)|\leq Cr\chi_{|\mu|}(t)$ for $|\mu|\leq 2\la-4$, $\mu_2+\mu_3=0$, and $r\leq\ds\frac{t}{2}+M+1$.
}

{\bf Proof.}
(1) If $t\leq 8M+7$, (1) can be got directly by Lemma 2.3 (4) and (6). Otherwise, we know that $t^2-r^2$ is equivalent to $t^2$ when $r\leq\ds\frac{3}{4}t+\frac{3}{2}(M+1)$. Set $m_\mu(t)=|\nabla\Gamma^\mu(\chi\dot\theta)(t)|+|\nabla\cdot\Gamma^\mu(\chi\dot w)(t)|$ and $M_\mu(t)=\ds\sum_{|\nu|\leq |\mu|}m_\nu(t)$. According to the expressions (2.7)-(2.9) and Lemma 2.3, we have
$$
m_\mu(t)\leq C\Bigg(\frac{\widetilde Q^-_{|\mu|+3}(t)}{\langle t\rangle^{3/2}}+\sum_{j=1}^{13}|\tilde f_j^\mu(t)|\bigg),
\eqno{(2.19)}
$$
here we have applied $|\tilde f_0^\mu(t)|+|\tilde f_{14}^\mu(t)|\leq C\langle t \rangle^{-3/2}E^{1/2}_{|\mu|+2}(t)\leq C\langle t \rangle^{-3/2}\widetilde Q^-_{|\mu|+3}(t)$.

We now treat $\ds\sum_{j=1}^{13}|\tilde f_j^\mu(t)|$.

If $\Gamma^\nu=\pa_t^lX^k$, then $\Gamma^\nu(\chi\dot w)(t,x)=
\ds\f{1}{r}\biggl(\int_0^r s(\nabla\cdot\Gamma^\nu(\chi\dot w))(s)ds\biggr)\frac{x}{r}$. At this time,
it follows the inequality
in Lemma 2.1 (c) of [1] that
$$
|\pa_x^\a\Gamma^\nu(\chi\dot w)(t)|\leq C_\a\sum_{j\leq |\a|-1}|\nabla\cdot\pa_x^j\Gamma^\nu(\chi\dot w)(t)|.\eqno{(2.20)}
$$

This, together with Lemma 2.3 and the expressions of $\tilde f_j^\mu$ ($1\le j\le 12$), yields
$$
\sum_{j=1}^{12}|\tilde f_j^\mu(t)|\leq C\Bigg(\ds\frac{\ve\widetilde Q^-_{|\mu|+2}(t)}{\langle t\rangle^3}+E_{|\mu|+2}^{1/2}(t)M_\mu(t)+\frac{\ve}{\langle t\rangle^{3/2}}M_\mu(t)\Bigg),\eqno{(2.21)}$$
or
$$
\sum_{j=1}^{12}|\tilde f_j^\mu(t)|\leq C\Bigg(\ds\frac{\ve\widetilde Q^-_{|\mu|+2}(t)}{\langle t\rangle^3}+\frac{\widetilde Q_{|\mu|+2}(t)}{\langle t\rangle^{1/2}}M_\mu(t)+\frac{\ve}{\langle t\rangle^{3/2}}M_\mu(t)\Bigg).\eqno{(2.22)}
$$

Noticing $\tilde f_{13}^\mu=-\pa^\mu (X+1)^k(\chi z\nabla\theta)+\pa^\a(X+1)^k(\chi\theta z\nabla\theta)$
and ${\widetilde Q}_{|\mu|+2}^-(t)\le C\ve$ for $|\mu|\le\la-2$ in terms of Lemma 2.4,
then by the assumption $\widetilde Q_{|\mu|+2}(t)\leq\eta\langle t\rangle^{1/2}$
for $|\mu|\geq\la-1$,  Lemma 2.3 and Lemma 2.2, we can obtain that $|\p^{\al'}X^{k'}\th(t)|$
is bounded for $|\al'|\le|\al|, k'\le k$ and further arrive at
\begin{align*}
|\tilde f_{13}^\mu(t)|\leq
&C\Big(\ve^2\langle t\rangle^{-5/2}+\ve M_\mu(t)+\ve\langle t\rangle^{-5/2}E^{1/2}_{|\mu|+2}(t)+\sum_{\nu\leq \mu}E^{1/2}_{|\nu|+2}(t)m_{\mu-\nu}(t)\\
&+\langle t\rangle^{-2}E_{[\frac{|\mu|}{2}]+2}^{1/2}\widetilde Q_{|\mu|+3}^-\Big).\tag{2.23}
\end{align*}

For $|\mu|\leq\la-2$, collecting (2.21), (2.23) yields
$$
\ds\sum_{j=1}^{13}|\tilde f_j^\mu(t)|\leq C\Big(\ds\f{\ve}{\langle t\rangle^{2}}{\widetilde Q}_{|\mu|+3}^-(t)
+\f{\ve^2}{\langle t\rangle^{5/2}}+\ve M_\mu(t)\Big),\eqno{(2.24)}
$$
which together with (2.19) means (1) holds.

For $\la-1\leq |\mu|\leq 2\la-4$, due to $E^{1/2}_{|\nu|+2}(t)m_{\mu-\nu}(t)\leq\ve M_\mu(t)$ if $|\nu|\leq\la-2$ and $E^{1/2}_{|\nu|+2}(t)m_{\mu-\nu}(t)\leq C\ve\langle t\rangle^{-3/2}E_{|\nu|+2}^{1/2}(t)$ otherwise, then by (2.22), (2.23)
$$
\sum_{j=1}^{13}|\tilde f_j^\mu(t)|\leq C\Big(\frac{\ve}{\langle t\rangle^2}\widetilde Q_{|\mu|+3}^-(t)+(\ve+\eta)M_\mu(t)+\frac{\ve^2}{\langle t\rangle^{5/2}}+\frac{\ve}{\langle t\rangle^{3/2}}E_{|\mu|+2}^{1/2}(t)\Big).
$$

This, together with (2.19), also yields (1).

(2) If $t$ is small, the estimate can be obtained by Lemma 2.3 easily.
Otherwise, (2) follows from (2.11) and (1).

(3) If $\mu_2+\mu_3\neq0$, as shown in (2.20), then we have for $|\mu|\leq 2\la-3$
$$
|\Gamma^\mu(\chi\dot w)(t)|\leq C\sum_{|\al|\leq \mu_2+\mu_3-1}|\nabla\cdot\pa_x^\al\pa_t^{\mu_1}X^{\mu_4}(\chi\dot w)(t)|\leq C\chi_{|\mu|-1}(t).
$$

(iv) If $\mu_2+\mu_3=0$, then $\Gamma^\mu\dot w(t,x)=\bigg(\ds\frac{1}{r}\int_0^rs(\nabla\cdot\pa_t^{\mu_1}X^{\mu_4}\dot w)(t,s){d s}\bigg)\ds\frac{x}{r}$ and (4) follows (1) immediately.  $\square$
\bigskip

When $\ve>0$ is small enough, as in [8], we next show that (1.3) is isentropic outside the fixed ball $B(0, M+1)$
for any time $t$.

{\bf Lemma 2.6.} {\it For small $\ve>0$, if $(\dot\theta, \dot w, \dot z)$ is a smooth solution to (2.4)
for $(t,x)\in [0, T]\times\Bbb R^2$ and $E_\la(t)\leq\ve^2$ with $\la\ge 4$, then $\dot z(t,x)=0$ for $|x|\geq M+1$.
}

{\it Proof.}
From the third equation in (2.3), we know that $\pa_tz+W(t,r)\pa_rz=0$, where $W(t,r)=w(t,r)\cdot\ds\frac{x}{r}$.
Define the characteristics $r=r(t)$ by $r'(t)=W(t, r(t))$ with $r(0)=M$, then it is easy to see that $z(t,x)=0$
for $|x|\geq r(t)$ by the compact support property of $z(0,x)$.

By Lemma 2.2 (1) and Lemma 2.3 (9) together with Lemma 2.4, we have $|w(t)|_-\leq |w_a(t)|_-+|\dot w(t)|_-
\le C\ve\langle t\rangle^{-3/2}+C\langle t\rangle^{-1/2}\widetilde Q_2^-(t)\leq C\ve\langle t\rangle^{-1/2}$.
In addition, it follows from Lemma 2.3 (2) that $|w(t, x)|\le |w_a(t,x)|+|\dot w(t,x)|\leq C\langle t\rangle^{-1/2}(\ve+E_2^{1/2}(t))\leq C\ve\langle t\rangle^{-1/2}$ holds for $|x|\geq \ds\frac{t}{2}+M+1$. Therefore,
$|r(t)|\leq (M+1)\langle t\rangle^{1/2}$. From this, we can take use of Lemma 2.2 (3) and Lemma 2.5 (4)
to get $|W(t, r(t))|\leq C\ve\langle t\rangle^{-1}$, and hence $|r(t)|\leq (M+1)\langle t\rangle^{\delta}$ for any $0<\delta<\frac{1}{2}$. Applying Lemma 2.2 (3) and Lemma 2.5 (4) again, we have $|r'(t)|\leq C\ve\langle t\rangle^{-3/2+\delta}$
and further $|r(t)-M|\leq C\ve$. Consequently, $z(t, r)\equiv0$ as $|x|\geq M+1$, so does $\dot z$.
\qquad \qquad $\square$\bigskip

\vskip 0.4 true cm \centerline{\bf $\S 3$. Analysis on $(\dot\th, \dot w, \dot z)$ near light cone
and establishment of a crucial energy estimate} \vskip 0.4 true cm

From Lemma 2.6, we know that $z(t,x)\equiv0$ and then $S\equiv\bar S$ holds when $|x|\geq M+1$.
Assume that $(\dot\theta, \dot w, \dot z)$ is a smooth solution to (2.4) for $0\leq t\leq T$, and let
$\phi\in C^\infty(\{(t,x)\in [0, T]\times\Bbb R^2: |x|\geq M+1\})$ be the corresponding potential
vanishing in $|x|\geq M+t$ and satisfying $\nabla\phi=u$, $\pa_t\phi=-\frac{1}{2}|u|^2-h(\rho)$.
For $(\lambda,y)$,
$(\tilde\lambda,\tilde y)\in\Bbb R\times\Bbb R^2$, we define
\begin{align*}
&F_1(\lambda,y)=\frac{1}{2}(\lambda^2-|y|^2),\\
&F_2((\lambda,y),(\tilde\lambda,\tilde y))=\frac{1}{2}(\lambda\tilde\lambda-y\cdot\tilde y).
\end{align*}

Let $\dot\phi=\phi-\phi_a$, $\xi=(\theta, w)$, $\xi_a=(\theta_a, w_a)$ and $\dot\xi=\xi-\xi_a=(\dot\theta, \dot w)$,
where $\phi_a$ and $(\theta_a, w_a)$  are given in (2.1) and Lemma 2.2 respectively. Then we have
\begin{align*}
&\theta_a(t,x)=-\pa_t\phi_a(t,x)+F_1(\xi_a)(t,x),\\
&\theta(t,x)=-\pa_t\phi(t,x)+F_1(\xi)(t,x)\qquad\text{for}\quad|x|\geq M+1,\\
&\dot\theta(t,x)=-\pa_t\dot\phi(t,x)+F_2(\dot\xi, 2\xi_a+\dot\xi)(t,x)\qquad\text{for}\quad|x|\geq M+1.
\end{align*}

In order to make use of both null conditions introduced in [4] for the second order quasilinear
wave equations, we will pay more attention to the forms of $F_1(\xi_a)$ and $F_2(\dot\xi, 2\xi_a+\dot\xi)$, that is, if $|x|\geq M+1$, then
\begin{align*}
&F_1(\xi_a)=\frac{1}{2}(\pa_t\phi_a)^2-\frac{1}{2}|\nabla\phi_a|^2+F_1(\xi_a)\big(\theta_a-\frac{1}{2}F_1(\xi_a)\big),\tag{3.1}\\
&F_2(\dot\xi, 2\xi_a+\dot\xi)=\frac{1}{2}\pa_t\dot\phi(2\pa_t\phi_a+\pa_t\dot\phi)-\frac{1}{2}\nabla\dot\phi
\cdot(2\nabla\phi_a+\nabla\dot\phi)+\frac{1}{2}F_2(\dot\xi,2\xi_a+\dot\xi)(2\theta_a+\dot\theta)\\
&\qquad \qquad \qquad\quad -\frac{1}{2}\big(F_2(\dot\xi,2\xi_a+\dot\xi)
-\dot\theta\big)\big(2F_1(\xi_a)+F_2(\dot\xi,2\xi_a+\dot\xi)\big).\tag{3.2}
\end{align*}

Next we cite two fundamental estimates on the first and second null conditions which are shown in [4], Lemma 6.64-Lemma 6.65
of [11] or [13]
(one can see Lemma 3.1-Lemma 3.5 of [13] for some details).

{\bf Lemma 3.1.} {\it
Assume that $g_{i}^{jk}\in\Bbb R$ with $g_i^{jk}=g_i^{kj}$ $(0\le i, j, k\le 2)$ and $\ds\sum_{i,j,k=0}^2g_i^{jk}p_ip_jp_k=0$ for any $p=(p_0,p_1,p_2)\in\Bbb R^3$ satisfying $p_0^2=p_1^2+p_2^2$. Then, there exists a positive constant $C$ such that for any $f, g\in C^\infty(
[0,T]\times\Bbb R^2)$,
$$
|\ds\sum_{i,j,k=0}^2g_i^{jk}(\pa_i f\pa_{jk}^2g)(t, x)|\leq C\big(|Zf(t,x)||\pa^2g(t,x)|
+|\pa f(t,x)||Z\pa g(t,x)|\big)\eqno{(3.3)}
$$
and
\begin{align*}
\sum_{|\al|\le n}|\Gamma^{\al}(\ds\sum_{i,j,k=0}^2g_i^{jk}\p_i f\p_{jk}^2f)(t,x)&
-\sum_{i,j,k=0}^2g_i^{jk}(\p_i f\p_{jk}^2\Gamma^{\al} f)(t,x)|\\
&\le C_n\ds\sum_{|\beta+\nu|\le n+1,\quad |\beta|,|\nu|\le n}|Z\Gamma^\beta f(t,x)||\Gamma^{\nu}\p f(t,x)|.\tag{3.4}\\
\end{align*}
}

{\bf Lemma 3.2.} {\it
Assume that $g_{ij}^{kl}\in\Bbb R$ with $g_{ij}^{kl}=g_{ij}^{lk}$ and $g_{ij}^{kl}=g_{ji}^{kl}$ ($0\le i, j, k, l\le 2$),
and $\ds\sum_{i,j,k,l=0}^2g_{ij}^{kl}p_ip_jp_kp_l$ $=0$ for any $p=(p_0,p_1,p_2)\in\Bbb R^3$ satisfying $p_0^2=p_1^2+p_2^2$.
Then, there exists a positive constant $C$ such that for any $f,g,h\in C^\infty([0, T]\times\Bbb R^2)$,
\begin{align*}
&|\ds\sum_{i,j,k,l=0}^2g_{ij}^{kl}(\pa_if\pa_jg\pa_{kl}^2h)(t,x)|\leq C\big(|\pa f(t,x)||\pa g(t,x)||Z\pa h(t,x)|\\
&\qquad +|\pa f(t,x)||Z g(t,x)||\pa^2h(t,x)|+|Zf(t,x)||\pa g(t,x)||\pa^2h(t,x)|\big).\tag{3.5}
\end{align*}
}

\bigskip

Based on Lemma 3.1 and Lemma 3.2, and noting $
|Z\varphi(t,x)|\leq C\ds\Big(\frac{||x|-t|}{t}|\pa\varphi(t,x)|+\frac{1}{t}(|X\varphi(t,x)|$ $+|\Omega\varphi(t,x)|)\Big)
$, then one can obtain the following estimates under the assumptions of Lemma 3.1 and Lemma 3.2:
\begin{align*}
&|\ds\sum_{i,j,k=0}^2g_i^{jk}(\pa_i\Gamma^{\al}\phi_a\pa_{jk}^2\Gamma^{\beta}\dot\phi)(t,x)|\le\frac{C\ve}{t\langle t\rangle^{1/2}}\sum_{|\nu|\le |\beta|+1}|\Gamma^{\nu}\pa\dot\phi(t,x)|,\tag{3.6}\\
&|\ds\sum_{i,j,k=0}^2g_i^{jk}(\pa_i\Gamma^\al\dot\phi\pa_{jk}^2\Gamma^\beta\phi_a)(t,x)|\le\frac{C\ve}{t\langle t\rangle^{1/2}}\bigg(\sum_{|\nu|\le |\al|+1}|(\si_-^{-1}\Gamma^\nu\dot\phi)(t,x)|
+\sum_{|\nu|\le |\al|}|\Gamma^\nu\pa\dot\phi(t,x)|\bigg),\tag{3.7}\\
&|\ds\sum_{i,j,k,l=0}^2g_{ij}^{kl}(\pa_i\Gamma^\a\phi_a\pa_j\Gamma^\beta\phi_a\pa_{kl}^2\Gamma^\mu \dot\phi)(t,x)|\le\frac{C\ve^2}{t\langle t\rangle}\sum_{|\nu|\le |\mu|+1}|\Gamma^\nu\pa\dot\phi(t,x)|,\tag{3.8}\\
&|\ds\sum_{i,j,k,l=0}^2g_{ij}^{kl}(\pa_i\Gamma^\a\phi_a\pa_j\Gamma^\beta \dot\phi\pa_{kl}^2\Gamma^\mu \phi_a)(t,x)|\le\frac{C\ve^2}{t\langle t\rangle}
\bigg(\ds\sum_{|\nu|\le |\beta|+1}|(\si_-^{-1}\Gamma^\nu\dot\phi)(t,x)|\\
&\qquad\qquad\qquad\qquad\qquad\qquad\qquad\qquad\qquad +\ds\sum_{|\nu|\le |\beta|}|\Gamma^\nu\pa\dot\phi(t,x)|\bigg)\tag{3.9}\\
\end{align*}
and
\begin{align*}
|\ds\sum_{i,j,k=0}^2&g_i^{jk}(\pa_i\Gamma^\al\dot\phi\pa_{jk}^2\Gamma^\beta\dot\phi)(t,x)|\le\frac{C}{t}\big(\sum_{|p|\le |\al|,\quad |q|\le|\beta|}|(\si_-\Gamma^{p}\pa\dot\phi\pa\Gamma^q\pa\dot\phi)(t,x)|\\
&+\sum_{|p|\le|\al|+1,\quad |q|\le |\beta|}|(\Gamma^p\dot\phi\pa\Gamma^q\pa \dot\phi)(t,x)|+\sum_{|p|\leq|\al|,\quad
|q|\le |\beta|+1}|(\Gamma^p \pa\dot\phi\Gamma^q\pa\dot\phi)(t,x)|\big),\tag{3.10}\\
|\ds\sum_{i,j,k,l=0}^2&g_{ij}^{kl}(\pa_i\Gamma^\a\dot\phi\pa_j\Gamma^\beta \phi_a\pa_{kl}^2\Gamma^\mu\dot\phi)(t,x)|\\
&\leq\frac{C\ve}{t\langle t\rangle^{1/2}}\big(\sum_{|p|\le |\a|,\quad |q|\le |\mu|+1}|(\Gamma^p\pa\dot\phi\Gamma^q\pa \dot\phi)(t,x)|
+\ds\sum_{ |p|\le |\a|+1\quad |q|\le |\mu|}|(\si_-^{-1}\Gamma^p\dot\phi\pa\Gamma^q\pa\dot\phi)(t,x)|\big),\tag{3.11}\\
\end{align*}
\begin{align*}
|\ds\sum_{i,j,k,l=0}^2&g_{ij}^{kl}(\pa_i\Gamma^\a\dot\phi\pa_j\Gamma^\beta\dot\phi\pa_{kl}^2\Gamma^\mu\phi_a)(t,x)|\\
&\leq\frac{C\ve}{t\langle t\rangle^{1/2}}\Big(\sum_{|p|\leq|\beta|+1,\quad |q|\leq|\a|}|(\si_-^{-1}\Gamma^p\dot\phi\Gamma^q\pa \dot\phi)(t,x)|+\ds\sum_{|p|\leq|\beta|,\quad |q|\leq|\a |+1}|(\si_-^{-1}\Gamma^q\dot\phi\Gamma^p\pa\dot\phi)(t,x)|\Big),\tag{3.12}\\
|\sum_{i,j,k,l=0}^2&g_{ij}^{kl}(\pa_i\Gamma^\a\dot\phi\pa_j\Gamma^\beta\dot\phi\pa_{kl}^2\Gamma^\mu\dot\phi)(t,x)|\\
&\leq\frac{C}{t}\Big(\sum_{|p|\leq|\a|,\quad |q|\leq|\beta|,|s|\leq|\mu|}|(\si_-\Gamma^p\pa\dot\phi\Gamma^q\pa \dot\phi\pa\Gamma^s\pa\dot\phi)(t,x)|\\
&+\sum_{|s|\leq|\mu|+1,\quad |q|\leq|\beta|,|p|\leq|\a|}|(\Gamma^p\pa\dot\phi\Gamma^q\pa\dot\phi\Gamma^s\pa\dot\phi)(t,x)|\\
\end{align*}

\begin{align*}
&+\sum_{|q|\leq|\beta|+1,\quad |s|\leq|\mu|,|p|\leq|\a|}|(\Gamma^p\pa \dot\phi\Gamma^q\dot\phi\pa\Gamma^s\pa\dot\phi)(t,x)|\\
&+\sum_{|p|\leq|\a|+1,\quad |q|\leq|\beta|,|s|\leq|\mu|}|(\Gamma^p\dot\phi\Gamma^q\pa\dot\phi\pa\Gamma^s\pa \dot\phi)(t,x)|\Big).\tag{3.13}
\end{align*}

To treat each term better in the right hand side of (3.6)-(3.13), we require the following two results
(Lemma 3.3 and Lemma 3.4):

{\bf Lemma 3.3.} {\it
Suppose that $f(t,x)\in C^1(\Bbb R_+\times\Bbb R^2)$ and $\text{supp} f\subset\{(t,x):|x|\leq M+t\}$.
Then there exists a positive constant $C$ independent of $t$ such that
$$
\|(\si_-^{-1}f)(t,\cdot)\|_{L^2(|x|\geq M+1)}\leq C\|\nabla f(t,\cdot)\|_{L^2(|x|\geq M+1)}.\eqno{(3.14)}
$$
}

{\bf Remark 3.1.} {\it Here we point out that the inequality $\|(\si_-^{-1}f)(t,\cdot)\|_{L^2}\leq C\|\nabla f(t,\cdot)\|_{L^2}$
in the whole space has been established in [15].}
\bigskip

{\bf Proof.}
Due to $\text{supp} f\subset\{|x|\leq M+t\}$, one then has
$$
f(t,x)=-\int_{|x|}^{M+t}\nabla f(t, s\omega)\cdot\omega  ds.
$$

This derives
$$
|f(t,x)|^2\leq\bigg(\int_{|x|}^{M+t}|\nabla f(t, s\omega)|^2(1+|t-s|)^{1/2} ds\bigg)
\int_{|x|}^{M+t}(1+|t-s|)^{-1/2} ds.
$$

With easy computations, we have $\int_{|x|}^{M+t}(1+|t-s|)^{-1/2} ds\leq C(1+t-|x|)^{1/2}$ for $t>|x|$, and $\int_{|x|}^{M+t}(1+|t-s|)^{-1/2} ds\leq C$ for $t\leq|x|$. Hence, if $t\geq M+1$, then
\begin{align*}
&\int_{M+1}^{M+t}|f(t, r\omega)|^2(1+|t-r|)^{-2}r dr\\
&\leq C\int_{M+1}^t\Big(\int_r^t|\nabla f(t, s\omega)|^2(1+t-s)^{1/2} ds\Big)(1+t-r)^{-3/2}r dr\\
&\quad +C\int_{M+1}^t\Big(\int_t^{M+t}|\nabla f(t, s\omega)|^2 ds\Big)(1+t-r)^{-3/2}r dr\\
&\quad +C\int_t^{M+t}\Big(\int_r^{M+t}|\nabla f(t, s\omega)|^2 ds\Big)r(1+r-t)^{-2} dr\\
&\leq C\int_{M+1}^{M+t}s|\nabla f(t, s\omega)|^2 ds.\tag{3.15}
\end{align*}

If $1\leq t\leq M+1$, then
\begin{align*}
&\int_{M+1}^{M+t}|f(t,r\omega)|^2(1+|t-r|)^{-2}r dr\leq C\int_{M+1}^{M+t}\int_r^{M+t}|\nabla f(t,s\omega)|^2r ds dr\\
&\leq C\int_{M+1}^{M+t}s|\nabla f(t, s\omega)|^2 ds.\tag{3.16}
\end{align*}

Combining (3.15) with (3.16) yields (3.14). \qquad\qquad
\qquad\qquad\qquad\qquad\qquad\qquad\qquad$\square$\bigskip

{\bf Lemma 3.4.} {\it
Suppose that $f(t,x)\in C^2(\Bbb R_+\times\Bbb R^2)$ and $\text{supp} f\subset\{(t, x):|x|\leq M+t\}$. Then there exists a positive constant $C$ independent of $t$ such that for $|x|\geq M+1$
$$
\si(x)^{1/2}\si_-(t,x)^{-1}|f(t,x)|\leq C\sum_{j=0}^1\sum_{|\al|=0}^{2-j}\|\pa\nabla^\al\Omega^jf(t,\cdot)\|_{L^2(|x|\geq M+1)}.\eqno{(3.17)}
$$
}
\bigskip

{\bf Proof.}
Similar to the proof of Lemma 1 (3.1) in [21], we have for $|x|\geq M+1$ and $g(t,x)\in C^2(\Bbb R_+\times\Bbb R^2)$
with $\text{supp} g\subset\{(t, x):|x|\leq M+t\}$
$$
\si(x)^{1/2}|g(t,x)|\leq C\sum_{j=0}^1\sum_{|\al|=0}^{2-j}\|\nabla^\al\Omega^jg(t,\cdot)\|_{L^2(|x|\geq M+1)}.\eqno{(3.18)}
$$

Choosing $g(t,x)=\si_-(t,x)^{-1}f(t,x)$ in (3.18) and applying Lemma 3.3 yield (3.17) directly.
$\square$\bigskip
Next, we extend Lemma 2.4 so that an energy estimate in the whole space is established as in Proposition 4.1 of [9].

{\bf Lemma 3.5.} {\it
For fixed $n\in\Bbb N$ with $n\geq 5$, if  $(\dot\theta, \dot w, \dot z)$ is a smooth solution of (2.4)
for $(t, x)\in [0, T]\times \Bbb R^2$, and $E_{[\frac{n+5}{2}]}^{1/2}(t)\leq\eta$ holds for small positive constant $\eta$,
then $\ds Q_n(t)\leq C(E_n^{1/2}(t)+\frac{\ve^2}{\langle t\rangle^{3/2}})$ for $0\leq t\leq T$.
}
\bigskip

{\bf Proof.} If $t\leq 1$, it is easy to see $Q_n(t)\leq C E_n^{1/2}(t)$. We now focus on the case of $t\geq 1$.

It is noted that for $|\nu|\leq n$
$$
\|\Gamma^\nu\pa \dot\phi(t)\|_+\leq C(\|\Gamma^\nu\dot\theta(t)\|_++\|\Gamma^\nu\dot w(t)\|_+)\leq C E^{1/2}_n(t),\eqno{(3.19)}
$$
here Lemma 2.3 (2) is applied in the first inequality.

Set $Q_n^+(t)=\ds\sum_{|\mu|\leq n-1}(\|\si_-(t)\nabla\Gamma^\mu\dot\theta(t)\|_++\|\si_-(t)\pa_t\Gamma^\mu\dot w(t)\|_++\|\si_-(t)\nabla\cdot\Gamma^\mu\dot w(t)\|_+)$. Similar to (2.12), the relationship between $E_j(t)$ and $Q_j^+(t)$
($1\le j\le n$) is
$$
Q_j^+(t)\leq C\big(E_j^{1/2}(t)+\sum_{|\mu|\leq j-1}t\|h_0^\mu(t)\|_++\sum_{|\mu|\leq j-1}\langle t\rangle\|h^\mu(t)\|_+\big).\eqno{(3.20)}
$$

When $|\mu|\leq j-1$, we have
$$
\|h_0^\mu(t)\|_+\leq\sum_{\nu\leq\mu}C_{\mu\nu}\big(\|I_1^{\mu\nu}(t)\|_++\|I_2^{\mu\nu}(t)\|_++\|I_3^{\mu\nu}(t)\|_+\big)
+\sum_{\nu\leq\mu}C_{\mu\nu}\|J_{\mu\nu}(t)\|_+,\eqno{(3.21)}
$$
where
\begin{align*}
I_1^{\mu\nu}(t,x)&=\sum_{i=1}^2\big((\Gamma^\nu\pa_i\phi_a)(\pa_i\Gamma^{\mu-\nu}\pa_t\dot\phi)
-(\Gamma^\nu\pa_t\phi_a)(\pa_i\Gamma^{\mu-\nu}\pa_i\dot\phi)\big)(t,x),\\
I_2^{\mu\nu}(t,x)&=\sum_{i=1}^2\big((\Gamma^\nu\pa_i\dot\phi)(\pa_i\Gamma^{\mu-\nu}\pa_t\phi_a)
-(\Gamma^\nu\pa_t\dot\phi)(\pa_i\Gamma^{\mu-\nu}\pa_i\phi_a)\big)(t,x),\\
I_3^{\mu\nu}(t,x)&=\sum_{i=1}^2\big((\Gamma^\nu\pa_i\dot\phi)(\pa_i\Gamma^{\mu-\nu}\pa_t\dot\phi)
-(\Gamma^\nu\pa_t\dot\phi)(\pa_i\Gamma^{\mu-\nu}\pa_i\dot\phi)\big)(t,x),
\end{align*}
and using (3.1)-(3.2) to get $J_{\mu\nu}(t,x)=\ds\sum_{j=1}^5\sum_{|\al|+|\beta|+|\gamma|\le\mu}C_{\al\beta\gamma}^{\nu}J^{\a\beta\gamma}_j(t,x)
+\tilde J_{\mu\nu}(t,x)$ with
\begin{align*}
J^{\a\beta\gamma}_1=\frac{1}{2}\sum_{k=1}^2&(\pa_t\Gamma^\a \phi_a)(\pa_t\Gamma^\beta \phi_a)(\pa_k^2\Gamma^\gamma \dot\phi)-\frac{1}{2}\sum_{k,l=1}^2(\pa_l\Gamma^\a \phi_a)(\pa_l\Gamma^\beta \phi_a)(\pa_k^2\Gamma^\gamma \dot\phi)\\
&-\sum_{k=1}^2(\pa_t\Gamma^\a \phi_a)(\pa_k\Gamma^\beta \phi_a)(\pa_k\pa_t\Gamma^\gamma \dot\phi)+\sum_{k,l=1}^2(\pa_l\Gamma^\a \phi_a)(\pa_k\Gamma^\beta \phi_a)(\pa_{kl}^2\Gamma^\gamma \dot\phi),
\end{align*}
\begin{align*}
J^{\a\beta\gamma}_2=-\sum_{k=1}^2&(\pa_t\Gamma^\beta \phi_a)(\pa_k\Gamma^\gamma\dot\phi)(\pa_k\pa_t\Gamma^\a \phi_a)+\sum_{k,l=1}^2(\pa_l\Gamma^\beta\phi_a)(\pa_k\Gamma^\gamma\dot\phi)(\pa_{kl}^2\Gamma^\a\phi_a)\\
&-\sum_{k=1}^2(\pa_k\Gamma^\beta\phi_a)(\pa_t\Gamma^\gamma\dot\phi)(\pa_k\pa_t\Gamma^\a\phi_a)+\sum_{k,l=1}^2(\pa_k\Gamma^\beta \phi_a)(\pa_l\Gamma^\gamma\dot\phi)(\pa_{kl}^2\Gamma^\a \phi_a)\\
&+\sum_{k=1}^2(\pa_t\Gamma^\beta\phi_a)(\pa_t\Gamma^\gamma\dot\phi)(\pa_k^2\Gamma^\a\phi_a)-\sum_{k,l=1}^2(\pa_k\Gamma^\beta \phi_a)(\pa_k\Gamma^\gamma\dot\phi)(\pa_l^2\Gamma^\a \phi_a),
\end{align*}
\begin{align*}
J^{\a\beta\gamma}_3=-\sum_{k=1}^2&(\pa_k\Gamma^\gamma\dot\phi)(\pa_t\Gamma^\beta\phi_a)(\pa_k\pa_t\Gamma^\a \dot\phi)+\sum_{k,l=1}^2(\pa_k\Gamma^\gamma\dot\phi)(\pa_l\Gamma^\beta\phi_a)(\pa_{kl}^2\Gamma^\a\dot\phi)\\
&-\sum_{k=1}^2(\pa_t\Gamma^\gamma\dot\phi)(\pa_k\Gamma^\beta\phi_a)(\pa_k\pa_t\Gamma^\a\dot\phi)+\sum_{k,l=1}^2(\pa_l \Gamma^\gamma \dot\phi)(\pa_k\Gamma^\beta\phi_a)(\pa_{kl}^2\Gamma^\a\dot\phi)\\
&+\sum_{k=1}^2(\pa_t\Gamma^\gamma\dot\phi)(\pa_t\Gamma^\beta\phi_a)(\pa_k^2\Gamma^\a \dot\phi)-\sum_{k,l=1}^2(\pa_k\Gamma^\gamma \dot\phi)(\pa_k\Gamma^\beta\phi_a)(\pa_l^2\Gamma^\a\dot\phi),
\end{align*}
\begin{align*}
J^{\a\beta\gamma}_4=-\sum_{k=1}^2&(\pa_k\Gamma^\gamma\dot\phi)(\pa_t\Gamma^\a\dot\phi)(\pa_k\pa_t\Gamma^\beta \phi_a)+\sum_{k,l=1}^2(\pa_k\Gamma^\gamma\dot\phi)(\pa_l\Gamma^\a\dot\phi)(\pa_{kl}^2\Gamma^\beta\phi_a)\\
&+\frac{1}{2}\sum_{k=1}^2(\pa_t\Gamma^\gamma\dot\phi)(\pa_t\Gamma^\a \dot\phi)(\pa_k^2\Gamma^\beta \phi_a)-\frac{1}{2}\sum_{k,l=1}^2(\pa_k\Gamma^\gamma\dot\phi)(\pa_k\Gamma^\a\dot\phi)(\pa_l^2\Gamma^\beta\phi_a),
\end{align*}
\begin{align*}
J^{\a\beta\gamma}_5=-\sum_{k=1}^2&(\pa_k\Gamma^\gamma\dot\phi)(\pa_t\Gamma^\beta\dot\phi)(\pa_k\pa_t\Gamma^\a \dot\phi)+\sum_{k,l=1}^2(\pa_k\Gamma^\gamma\dot\phi)(\pa_l\Gamma^\beta\dot\phi)(\pa_{kl}^2\Gamma^\a\dot\phi)\\
&+\frac{1}{2}\sum_{k=1}^2(\pa_t\Gamma^\gamma\dot\phi)(\pa_t\Gamma^\beta\dot\phi)(\pa_k^2\Gamma^\a \dot\phi)-\frac{1}{2}\sum_{k,l=1}^2(\pa_k\Gamma^\gamma\dot\phi)(\pa_k\Gamma^\beta\dot\phi)(\pa_l^2\Gamma^\a\dot\phi),
\end{align*}
and $\tilde J_{\mu\nu}$ is a higher order error term whose explicit expression is not needed.

In terms of Lemma 3.3 and (3.19), it is easy to find that
\begin{align*}
\|I_1^{\mu\nu}(t)\|_++\|I_2^{\mu\nu}(t)\|_+&\leq\frac{C\ve}{\langle t\rangle^{3/2}}E_j^{1/2}(t),\tag{3.22}\\
\|J_1^{\a\beta\gamma}(t)\|_++\|J_2^{\a\beta\gamma}(t)\|_+&\leq\frac{C\ve^2}{\langle t\rangle^2}E_j^{1/2}(t).\tag{3.23}
\end{align*}

To deal with $I_3^{\mu\nu}$, we should pay attention to (3.10). If $|\al|\leq |\beta|$ and $|\al+\beta|\leq j-1$,
then we can obtain with the help of Lemma 3.4
\begin{align*}
&\sum_{|p|\leq|\al|+1,\quad |q|\leq|\beta|}\|(\Gamma^p\dot\phi)(\pa\Gamma^q\pa\dot\phi)(t)\|_+\leq\sum_{ |p|\leq|\al|+1,\quad |q|\leq|\beta|}|\si_-^{-1}\Gamma^p\dot\phi(t)|_+\|\si_-\pa\Gamma^q\pa\dot\phi(t)\|_+\\
&\leq\frac{C}{\langle t\rangle^{1/2}}E_{[\frac{j+5}{2}]}^{1/2}(t)\widetilde Q_j(t).
\end{align*}
If $|\al|>|\beta|$, by using Lemma 2.3 and Lemma 3.3 we have
\begin{align*}
&\sum_{|p|\leq|\al|+1,\quad |q|\leq|\beta|}\|(\Gamma^p \dot\phi)(\pa\Gamma^q\pa\dot\phi)(t)\|_+\leq\sum_{ |p|\leq|\al|+1,\quad |q|\leq|\beta|}|\si_-\pa\Gamma^q\pa \dot\phi(t)|_+\|\si_-^{-1}\Gamma^p \dot\phi(t)\|_+\\
&\leq\frac{C}{\langle t\rangle^{1/2}}\widetilde Q_{[\frac{j}{2}]+2}(t)E_{j}^{1/2}(t).
\end{align*}

Similarly, we can arrive at for $|\al+\beta|\leq j-1$
\begin{align*}
\sum_{|p|\leq|\al|,\quad |q|\leq|\beta|}\|\si_-(\Gamma^{p}\pa \dot\phi)(\pa\Gamma^q\pa \dot\phi)(t)\|_+&+\sum_{ |p|\leq|\al|,\quad |q|\leq|\beta|+1}\|(\Gamma^p \pa\dot\phi)(\Gamma^q\pa \dot\phi)(t)\|_+\\
\leq\frac{C}{\langle t\rangle^{1/2}}&\Big(E^{1/2}_{[\frac{j+3}{2}]}\widetilde Q_j(t)+\widetilde Q_{[\frac{j}{2}]+2}(t)E_{j-1}^{1/2}(t)\Big),
\end{align*}
$$
\sum_{|p|\leq |\beta|+1,\quad |q|\leq |\al|}\|\si_-^{-1}\Gamma^p\dot\phi\Gamma^q\pa \dot\phi(t)\|_+\leq\frac{C}{\langle t\rangle^{1/2}}\Big(E^{1/2}_{[\frac{j}{2}]+2}(t)E^{1/2}_{j-1}(t)+E^{1/2}_{[\frac{j+3}{2}]}(t)E_j^{1/2}(t)\Big).
$$

Therefore,
\begin{align*}
\|I_3^{\mu\nu}(t)\|_+&\leq\frac{C}{\langle t\rangle^{3/2}}\Big(E^{1/2}_{[\frac{j+5}{2}]}(t)\widetilde Q_j(t)+\widetilde Q_{[\frac{j}{2}]+2}(t)E_j^{1/2}(t)\Big),\tag{3.24}\\
\|J_3^{\a\beta\gamma}(t)\|_+&\leq\frac{C\ve}{\langle t\rangle^2}\Big(E^{1/2}_{[\frac{j+3}{2}]}(t)\widetilde Q_j(t)+\widetilde Q_{[\frac{j}{2}]+2}(t)E_{j-1}^{1/2}(t)\Big),\tag{3.25}\\
\|J_4^{\a\beta\gamma}(t)\|_+&\leq\frac{C\ve}{\langle t\rangle^2}\Big(E^{1/2}_{[\frac{j}{2}]+2}(t)E^{1/2}_{j-1}(t)+ E^{1/2}_{[\frac{j+3}{2}]}(t)E_j^{1/2}(t)\Big).\tag{3.26}
\end{align*}

Analogously, we can also get the estimate of $J_5^{\a\beta\gamma}$ by means of
comparing the size of $|p|$, $|q|$ and $|s|$ in (3.13), that is,
\begin{align*}
&\|J_5^{\a\beta\gamma}(t)\|_+\\
&\leq\frac{C}{\langle t\rangle^2}\Big(\widetilde Q_{[\frac{j+5}{2}]}(t)E^{1/2}_{[\frac{j+3}{2}]}(t)E_{j-1}^{1/2}(t)+E_{[\frac{j}{2}]+1}(t)\widetilde Q_j(t)+E_{[\frac{j+5}{2}]}E^{1/2}_{j-1}(t)+E_{[\frac{j}{2}]+1}(t)E_j^{1/2}(t)\\
&\qquad +E_{[\frac{j+5}{2}]}^{1/2}(t)\widetilde Q_{[\frac{j+5}{2}]}(t)E^{1/2}_{j-1}(t)+E_{[\frac{j+3}{2}]}^{1/2}(t)\widetilde Q_{[\frac{j+5}{2}]}(t)E^{1/2}_{j}(t)+E^{1/2}_{[\frac{j}{2}]+1}(t)E^{1/2}_{[\frac{j}{2}]+2}(t)\widetilde Q_j(t)\Big).\tag{3.27}
\end{align*}

In addition, it is noted that $\tilde J_{\mu\nu}$ only contains the high-order error terms, then by a direct verification
one can derive that the $L^2$ norm of $\tilde J_{\mu\nu}$ near the light cone can be controlled by those terms in the right hand sides of (3.22)-(3.27).

Due to $Q_n^-(t)\leq C\big(E_n^{1/2}(t)+\ds\frac{\ve^2}{\langle t\rangle^{3/2}}\big)$ in terms of Lemma 2.4,
we apply the estimates (3.21)-(3.27) to obtain  $\|h_0^\mu(t)\|_+\leq \ds\frac{C}{\langle t\rangle^{3/2}}\eta(\eta+\frac{\ve^2}{\langle t\rangle^{3/2}})+\frac{C\eta}{\langle t\rangle^{3/2}}Q^+_{[\frac{n+5}{2}]}(t)$ when $|\mu|\leq [\ds\frac{n+5}{2}]-1$ and $\eta$ is small enough. In addition, because of $z(t,x)\equiv0$ for $|x|\geq M+1$, then we can get the same result for $h^\mu$ by the analogous analysis to $h_0^{\mu}$, that is, $\|h^\mu(t)\|_+\leq \ds\frac{C}{\langle t\rangle^{3/2}}\eta(\eta+\frac{\ve^2}{\langle t\rangle^{3/2}})+\frac{C\eta}{\langle t\rangle^{3/2}}Q^+_{[\frac{n+5}{2}]}(t)$ for $|\mu|\leq [\ds\frac{n+5}{2}]-1$. Hence,
$Q^+_{[\frac{n+5}{2}]}(t)\leq C\eta$ can be derived by utilizing (3.20).
Taking $j=n$ in the estimates in (3.22)-(3.27) and applying (3.20) again, we obtain
$$
Q^+_{n}(t)\leq C\big(E_n^{1/2}(t)+\frac{\ve^2}{\langle t\rangle^{5/2}}\big).
$$

This, together with Lemma 2.4, yields Lemma 3.5.
\qquad \qquad $\square$
\bigskip

\vskip 0.4 true cm \centerline{\bf $\S 4$. Global energy estimates and proof of Theorem 1.1} \vskip 0.4 true cm

From (2.4), as in [9], we have
\begin{align*}
&(\pa_t+w\cdot\nabla)\Gamma^\mu\dot\theta+(1-\theta)\nabla\cdot\Gamma^\mu\dot w=\hat h_0^\mu\equiv
\ds\sum_{j\in\{2,5\}}f_j^\mu+\sum_{j\in\{1,3,4,6\}}\hat f_j^\mu,\tag{4.1}\\
&\ds(\frac{1}{1+z}\pa_t+\frac{w}{1+z}\cdot\nabla)\Gamma^\mu\dot w+(1-\theta)\nabla\Gamma^\mu\dot\theta=
\frac{\hat h^\mu}{1+z}\equiv\ds\f{1}{1+z}\bigg(\sum_{j\in\{8,11\}}f_j^\mu+\sum_{j\in\{7,9,10,12,13\}}\hat f_j^\mu\bigg),\tag{4.2}\\
&(\pa_t+w\cdot\nabla)\Gamma^\mu\dot z=\hat g^\mu,\tag{4.3}
\end{align*}
where $f_j^{\mu} (j=2, 5, 8, 11)$ have been defined in (2.6); if $j\not=13$, then $\hat f_j^0=0$, and $\hat f_j^\mu$ with $\mu\not=0$ are defined as $f_j^\mu$  but with the supplementary restriction condition $\nu\not=0$ in the sum; if $j=13$, then
$\hat f_j^\mu=f_j^\mu+(1-\theta)z\nabla\Gamma^\mu\dot\theta$. In addition, $\hat g^\mu=\ds\sum_{j\in\{1,2\}}\hat g_{aj}^\mu+\ds\sum_{j\in\{1,2\}}\hat g_{j}^\mu$ with $\hat g_{a1}^\mu=-\ds\sum_{\nu\leq \mu}\binom{\mu}{\nu}\Gamma^\nu w_a\cdot\nabla\Gamma^{\mu-\nu}z(0,x)$,
$ \hat g_{a2}^\mu =-\ds\sum_{0<\nu\leq \mu}\binom{\mu}{\nu}\Gamma^\nu w_a\cdot\nabla\Gamma^{\mu-\nu}\dot z$
and $\hat g_{1}^\mu=-\ds\sum_{\nu\leq \mu}\binom{\mu}{\nu}\Gamma^\nu \dot w\cdot\nabla\Gamma^{\mu-\nu}z(0,x)$,
$\hat g_{2}^\mu =-\ds\sum_{0<\nu\leq \mu}\binom{\mu}{\nu}\Gamma^\nu \dot w\cdot\nabla\Gamma^{\mu-\nu}\dot z$
if $\mu\neq0$, otherwise, $\hat g_{a2}^\mu= \hat g_{2}^\mu=0$ for $\mu=0$.

Set

$$\zeta^\mu=
\left(
\begin{array}{ccc}
\Gamma^\mu\dot\theta\\
\Gamma^\mu\dot w\\
\Gamma^\mu\dot z
\end{array}
\right), \quad
F^\mu= \left(
\begin{array}{ccc}
\hat h_0^\mu\\
\ds\frac{\hat h^\mu}{1+z}\\
\hat g^\mu
\end{array}
\right), \quad
A_0= \left(
\begin{array}{cccc} 1 & 0 & 0 & 0 \\ 0
& \ds\frac{1}{1+z} & 0 & 0 \\
0  & 0  & \ds\frac{1}{1+z} & 0  \\
0 & 0 & 0 & 1
\end{array}
\right),
$$

$$
A_1= \left(
\begin{array}{cccc}
w_1 & 1-\theta & 0 & 0 \\
1-\theta & \ds\frac{w_1}{1+z} & 0 & 0 \\
0  & 0  & \ds\frac{w_1}{1+z} & 0  \\
0 & 0 & 0 & w_1
\end{array}
\right),\quad
A_2= \left(
\begin{array}{cccc}
w_2 & 0 & 1-\theta & 0 \\
0 & \ds\frac{w_2}{1+z} & 0 & 0 \\
1-\theta  & 0  &\ds \frac{w_2}{1+z} & 0  \\
0 & 0 & 0 & w_2
\end{array}
\right),
$$
where $w=(w_1, w_2)$.

In this case, (4.1)-(4.3) can be written as $A_0\pa_t\zeta^\mu+A_1\pa_1\zeta^\mu+A_2\pa_2\zeta^\mu=F^\mu$.
Taking inner product with $\zeta^\mu$ in the space $\Bbb R^2$ and applying integration by parts, we arrive at
$$
\frac{ d}{d t}\langle A_0\zeta^\mu,\zeta^\mu\rangle(t)=2\langle F^\mu,\zeta^\mu\rangle(t)+\sum_{0\leq j\leq2}\langle(\pa_jA_j)\zeta^\mu,\zeta^\mu\rangle(t).\eqno{(4.4)}
$$

As in [9], we set $H_j^\mu=\langle f_j^\mu,\Gamma^\mu\dot\theta\rangle$ if $j=2,5$, $H_j^\mu=\langle\hat f_j^\mu,\Gamma^\mu\dot\theta\rangle$ if $j=1,3,4,6$,
$H_j^\mu=\langle\ds\frac{f_j^\mu}{1+z},\Gamma^\mu\dot w\rangle$ if $j=8,11$, $H_j^\mu=\langle\ds\frac{\hat f_j^\mu}{1+z},\Gamma^\mu\dot w\rangle$ if $j=7,9,10,12,13$,
$\hat H_{aj}^\mu=\langle\hat g_{aj}^\mu,\Gamma^\mu\dot z\rangle$ and $\hat H_{j}^\mu
=\langle\hat g_{j}^\mu,\Gamma^\mu\dot z\rangle$ for $j=1,2$. Then
$$
\langle F^\mu,\zeta^\mu\rangle=\sum_{1\leq j\leq 13}H_j^\mu+\sum_{j=1,2}\hat H_{aj}^\mu
+\sum_{j=1,2}\hat H_{j}^\mu.\eqno{(4.5)}
$$

Although the local existence of solution to (2.4) has  been early established (for example, one can see [18]),
we still give some detailed illustrations for later uses.

{\bf Lemma 4.1.} {\it
If $\ve>0$ is small, then (2.4) has a unique smooth solution $(\dot\theta, \dot w, \dot z)$ for $t\leq2/\ve$ and $E_m^{1/2}(t)\leq C\ve^2$.
}

\bigskip
{\bf Proof.} Define
$$
T_1=\sup\{T>0: E_m^{1/2}(t)\leq \ve,0<t\leq 2/\ve\}\quad\text{for}\quad m\geq5.
$$

It is noted that $|z(t,x)|\le C\ve$ holds for $t\in [0, T_1]$ by Lemma 2.3 (2). This, together with
Lemma 2.3 and Lemma 3.5, yields
$$
|\sum_{j=0}^2\pa_jA_j(t,x)|\leq C\big(|\nabla\theta(t,x)|+|\nabla\cdot w(t,x)|
+|w\cdot\nabla z(t,x)|\big)\leq\frac{C}{\langle t\rangle^{1/2}}(\ve+\widetilde Q_3(t))\leq\frac{ C\ve}{\langle t\rangle^{1/2}}.
$$

In addition, a direct computation yields
$$
|H_j^\mu(t)|\leq\frac{C\ve}{\langle t\rangle^{1/2}}E_m(t)\quad\text{for}\quad|\mu|\leq m\quad\text{and}\quad j\in\{1,2,4,5,7,8,10,11\}.\eqno{(4.6)}
$$

With respect to $H_3^\mu(t)$ for $|\mu|\leq m$, we will treat it under two kinds of cases:

If $|\nu|\leq |\mu-\nu|\leq m-1$, then by Lemma 2.3 and Lemma 3.5
$$
|\langle\Gamma^\nu\dot w\cdot\nabla\Gamma^{\mu-\nu}\dot\theta,\Gamma^\mu\dot\theta\rangle(t)|
\leq\frac{C}{\langle t\rangle^{1/2}}|\langle\si^{1/2}\Gamma^\nu\dot w\cdot\si_-\nabla\Gamma^{\mu-\nu}\dot\theta,\Gamma^\mu\dot\theta\rangle(t)|\le \ds\frac{C\ve}{\langle t\rangle^{1/2}}E_m(t)+\frac{C\ve^3}{\langle t\rangle^2}E_m^{1/2}(t).
$$

Similarly, for $|\nu|> |\mu-\nu|$,
$$
|\langle\Gamma^\nu\dot w\cdot\nabla\Gamma^{\mu-\nu}\dot\theta,\Gamma^\mu\dot\theta\rangle(t)|
\leq\frac{C}{\langle t\rangle^{1/2}}|\langle\si^{1/2}\si_-\nabla\Gamma^{\mu-\nu}\dot\theta\cdot\Gamma^\nu\dot w,\Gamma^\mu\dot\theta\rangle(t)|\leq\frac{C\ve}{\langle t\rangle^{1/2}}E_m(t).
$$

Therefore,
$$
|H_3^\mu(t)|\leq\frac{C\ve}{\langle t\rangle^{1/2}}E_m(t)+\frac{C\ve^3}{\langle t\rangle^2}E_m^{1/2}(t).\eqno{(4.7)}
$$

By the same way, we have
$$
|H_j^\mu(t)|\leq\frac{C\ve}{\langle t\rangle^{1/2}}E_m(t)+\frac{C\ve^3}{\langle t\rangle^2}E_m^{1/2}(t)
\quad\text{for}\quad|\mu|\leq m\quad\text{and}\quad j\in\{6,9,12\}.\eqno{(4.8)}
$$

On the other hand,
$$
\|\pa^\a(X+1)^k\big((1-\theta)z(0,x)\nabla\theta_a\big)\|
+\|\pa^\a(X+1)^k\big((1-\theta)\dot z\nabla\theta_a\big)\|\leq\frac{C\ve^2}{\langle t\rangle^{5/2}}+\frac{C\ve}{\langle t\rangle^{5/2}}E^{1/2}_m(t)\eqno{(4.9)}
$$
and
\begin{align*}
\|\pa^\a(X+1)^k\big((1-\theta)z\nabla\dot\theta\big)-(1-\theta)z\nabla\Gamma^\mu\dot\theta\|&\leq \frac{C}{\langle t\rangle}\sum_{0<\nu\leq \mu}\|\Gamma^\nu\big((1-\theta)z\big)\cdot\si_-\nabla\Gamma^{\mu-\nu}\dot\theta\|\\
&\le \ds\frac{C\ve}{\langle t\rangle}E^{1/2}_m(t)
+\frac{C\ve^3}{\langle t\rangle^{5/2}},\tag{4.10}
\end{align*}
which mean
$$
|H_{13}^\mu(t)|\leq\frac{C\ve^2}{\langle t\rangle^{5/2}}E^{1/2}_m(t)+\frac{C\ve}{\langle t\rangle}E_m(t).\eqno{(4.11)}
$$

Moreover, one can easily obtain
$$\|\hat g_{a1}^\mu(t)\|+\|\hat g_{1}^\mu(t)\|\leq\frac{C\ve^2}{\langle t\rangle^{5/2}}+ C\ve E_m^{1/2}(t),
$$
$$
\|\hat g_{a2}^\mu(t)\|+\|\hat g_{2}^\mu(t)\|\leq\frac{C\ve}{\langle t\rangle^{5/2}}E_m^{1/2}(t)+C\ve E_m^{1/2}(t)
$$
and then
$$
\sum_{j=1,2}|\hat H_{aj}^\mu|+\sum_{j=1,2}|\hat H_{j}^\mu|\leq\frac{C\ve^2}{\langle t\rangle^{5/2}}E_m^{1/2}(t)
+C\ve E_m(t).\eqno{(4.12)}
$$

By (4.6)-(4.9) and (4.11)-(4.12), we have from  (4.4)
$$
\frac{d}{dt}E_m(t)\leq\frac{C\ve^2}{\langle t\rangle^{2}}E_m^{1/2}(t)+C\ve E_m(t),
$$
and then for sufficiently small $\ve$,
$$
E_m^{1/2}(t)\leq C\ve^2e^{C\ve t}\leq C\ve^2.\qquad\qquad\qquad\qquad \qquad\qquad
$$

Therefore, Lemma 4.1 is proved by the local existence of solution and continuous induction method.$\square$

\bigskip

As in Lemma 4.1, in order to prove Theorem 1.1 by the continuous induction method, we require to
establish a uniform estimate on the solution $(\rho, w, z)$. To this end, we will derive the a
priori estimate on the related potential $\phi$ in the domain
$\{x: |x|\ge M+1\}$.

{\bf Lemma 4.2.} {\it
Assume that $k$ and $\la$ are integers with $[\ds\frac{k+1}{2}]+3\leq\la\leq k$ and $k\geq7$,
$(\dot\theta, \dot w, \dot z)$ is a smooth solution of (2.4) for $(t,x)\in [\ds\f{1}{\ve},  T]\times\Bbb R^2$.
If $E_\la(t)\leq\ve^2$, then one can find a positive number $C$ independent of $\ve$ and $T$
such that the potential $\phi(t,x)$ of velocity $w(t,x)$ in the domain $\{x: |x|\ge M+1\}$ satisfies
\begin{align*}
&\sum_{|\mu|\leq k}\int_{1/\ve}^t\int_{D_+(\tau)}\si_-(t,x)^{-2}|Z\Gamma^\mu \phi(\tau, x)|^2d\tau dx \\
&\qquad \leq C\Big(\ve^4+\sum_{|\mu|\leq k}\int_{1/\ve}^t\langle\tau\rangle^{-1}\int_{\{|x|\geq M+1\}}|\Gamma^\mu\pa\phi(\tau, x)|^2d\tau dx \Big)\tag{4.13}
\end{align*}
}\bigskip

{\bf Remark 4.1.} {\it In order to apply energy integral method to derive (4.13), we require to choose a
different ``ghost weight'' from the one in [4] due to the following reason: notice that the both null conditions
of 2-D quasilinear wave equation are fulfilled in the whole space $\Bbb R^2$ in [4], but for our 2-D compressible
Euler system (1.3),
the null conditions  hold only in the exterior domain $\{|x|\ge M+1\}$, and thus it is natural for us to multiply
a smooth cut-off function on the potential function $\phi$ so that a suitable weighted energy estimate can be obtained.
Due to this way of doing, the resulted ``ghost weight'' should be
reconsidered by comparison with that in [4]. More concretely speaking, the author in [4]
can obtain the a priori estimate $|\p\G^\mu v(t,x)|\le C\ve\si_-^{-1}(t,x)\si^{-\f12}(t,x)$ for the solution $v$ of
quasilinear wave equation $\p_t^2v-\Delta v+\ds\sum_{0\le i, j\le 2}g_{ij}(\p v)\p_{ij}^2v=0$ when the both null conditions hold, however, here we only get $|\p\G^\mu \phi(t,x)|\le C\ve\si_-^{-\f12}(t,x)\si^{-\f12}(t,x)$ for the potential $\phi$
(see (4.17) below) which will lead to a different choice of the ghost weight from that in [4].}

\bigskip

{\bf Proof.}
As in (2.1),  the function $\phi(t,x)$ satisfies for $|x|\geq M+1$
$$
\displaystyle\pa_t^2\phi-\triangle \phi+2\sum_{j=1}^2\pa_j\phi\pa_t\pa_j\phi+\sum_{j,k=1}^2\pa_j\phi\pa_k\phi\pa_{jk}^2\phi
-(2\pa_t\phi+|\nabla\phi|^2)\triangle\phi=0.\eqno{(4.14)}
$$

Define
$$
\varphi(t,x)=\vp(|x|-t)=\int_{-\infty}^{|x|-t}\frac{1}{(1+|\rho|)^{3/2}} d\rho,
$$
then $0\leq\varphi\leq4$ and $\vp'(|x|-t)=\si_-(t,x)^{-3/2}$ hold. In addition, we set $\tilde\si(t,x)=\big(1+(|x|-t-M)^2\big)^{1/2}$ and $\chi(t,x)=\tilde\chi(\ds\frac{2r}{t+2M+2})$ with the smooth function
\begin{equation*}
\tilde\chi(s) =
\left\{
\begin{aligned}
0, & \quad \text{if $s\leq\frac{1}{2}$},\\
1, & \quad \text{if $s\geq1$}.
\end{aligned}
\right.
\end{equation*}
Let the corresponding energy be denoted by
$$
\widetilde E_n(t)=\ds\sum_{|\mu|\le n}\int_{\Bbb R^2}\tilde\si(t,x)^{-1/2}e^{\varphi(t,x)}\chi(t,x)\big(|\pa_t\G^\mu \phi(t,x)|^2+|\nabla\G^\mu\phi(t,x)|^2\big) dx
$$
for $n\in\Bbb N\cup\{0\}$. Motivated by the terminology in [4], the weight function $\t\si(t,x)^{-1/2}e^{\varphi(t,x)}\chi(t,x)$
is also called the ghost weight by us, which will display the null conditions and decay rate simultaneously.

Notice that
\begin{align*}
\ds\pa_t^2\Gamma^\mu \phi-&\triangle\Gamma^\mu  \phi+2\sum_{j=1}^2\pa_j\phi\pa_t\pa_j\Gamma^\mu \phi+\sum_{j,k=1}^2\pa_j\phi\pa_k\phi\pa_{jk}^2\Gamma^\mu \phi-(2\pa_t\phi+|\nabla\phi|^2)\triangle \Gamma^\mu \phi\\
&=\sum_{\nu\leq \mu}C_\nu\Gamma^\nu\big(-2\sum_{j=1}^2\pa_j\phi\pa_t\pa_j\phi-\sum_{j,k=1}^2\pa_j\phi\pa_k\phi\pa_{jk}^2\phi
+(2\pa_t\phi+|\nabla\phi|^2)\triangle\phi\big)\\
&+2\sum_{j=1}^2\pa_j\phi\pa_t\pa_j\Gamma^\mu \phi+\sum_{j,k=1}^2\pa_j\phi\pa_k\phi\pa_{jk}^2\Gamma^\mu
\phi-(2\pa_t\phi+|\nabla\phi|^2)\triangle \Gamma^\mu \phi.\tag{4.15}
\end{align*}

Multiplying (4.15) by $\tilde\si^{-1/2}(t,x)e^{\varphi(t,x)}\chi(t,x)\pa_t\Gamma^\mu \phi$, integrating in the space $\Bbb R^2$ and using integration by parts, we can get
\begin{align*}
&\frac{1}{2}\frac{d}{dt}\int_{\Bbb R^2}\tilde\si^{-1/2}e^{\varphi}\chi\big(|\pa_t\Gamma^\mu \phi|^2+|\nabla\Gamma^\mu \phi|^2\big) dx+\frac{1}{2}\int_{\Bbb R^2}\tilde\si^{-1/2}e^{\varphi}\chi\varphi'|Z\Gamma^\mu \phi|^2 dx\\
&+\frac{ d}{dt}\int_{\Bbb R^2}\tilde\si^{-1/2}e^{\varphi}\chi\Bigg(\pa_t\phi|\nabla\Gamma^\mu\phi|^2-\frac{1}{2}\sum_{i,j=1}^2\pa_i \phi\pa_j\phi(\pa_j\Gamma^\mu\phi)(\pa_i\Gamma^\mu\phi)+\frac{1}{2}|\nabla\phi|^2|\nabla\Gamma^\mu\phi|^2\Bigg)dx\\
&-\frac{1}{4}\int_{\Bbb R^2}(|x|-t-M)\tilde\si^{-5/2}e^{\varphi}\chi|Z\Gamma^\mu\phi|^2 dx+\sum_{i=2}^5L^\mu_i(t)
=L^\mu_1(t),\tag{4.16}
\end{align*}

where
\begin{align*}
L^\mu_1(t)=&\int_{\Bbb R^2}\tilde\si^{-1/2}e^{\varphi}\chi\pa_t\Gamma^\mu\phi\Big\{\sum_{\nu\leq \mu}C_{\nu}\Gamma^\nu\bigg(-2\sum_{j=1}^2\pa_j\phi\pa_t\pa_j\phi-\sum_{j,k=1}^2\pa_j\phi\pa_k\phi\pa_{jk}^2\phi+(2\pa_t\phi\\
&+|\nabla\phi|^2)\triangle \phi\bigg)+2\sum_{j=1}^2\pa_j\phi\pa_t\pa_j\Gamma^\mu\phi+\sum_{j,k=1}^2\pa_j\phi\pa_k\phi\pa_{jk}^2\Gamma^\mu\phi
-(2\pa_t\phi+|\nabla\phi|^2)\triangle \Gamma^\mu\phi\Big\} dx,\\
L^\mu_2(t)=&-\int_{\Bbb R^2}\tilde\si^{-1/2}e^{\varphi}\varphi'\chi\Big\{\sum_{i=1}^2\omega_i\pa_i \phi(\pa_t\Gamma^\mu\phi)^2-2\sum_{i=1}^2\omega_i\pa_t\phi(\pa_i\Gamma^\mu\phi)(\pa_t\Gamma^\mu\phi)\\
&-\sum_{i=1}^2\pa_t\phi(\pa_i\Gamma^\mu\phi)^2+\sum_{i,j=1}^2\omega_i\pa_i\phi\pa_j\phi(\pa_j\Gamma^\mu\phi)
(\pa_t\Gamma^\mu\phi)+\frac{1}{2}\sum_{i,j=1}^2\pa_i\phi\pa_j\phi(\pa_i\Gamma^\mu\phi)(\pa_j\Gamma^\mu\phi)\\
&-\sum_{i=1}^2\omega_i|\nabla\phi|^2(\pa_i\Gamma^\mu\phi)(\pa_t\Gamma^\mu\phi)-\frac{1}{2}|\nabla\phi|^2|\nabla\Gamma^\mu \phi|^2\Big\}dx,\\
\end{align*}
\begin{align*}
L^\mu_3(t)=&-\int_{\Bbb R^2}\tilde\si^{-1/2}e^{\varphi}\chi\Big\{\sum_{i=1}^2\pa_i^2\phi(\pa_t\Gamma^\mu\phi)^2
-2\sum_{i=1}^2\pa_i\pa_t\phi(\pa_i\Gamma^\mu \phi)(\pa_t\Gamma^\mu\phi)+\pa_t^2\phi|\nabla\Gamma^\mu\phi|^2\\
&+\sum_{i,j=1}^2\pa_i(\pa_i\phi\pa_j\phi)(\pa_j\Gamma^\mu\phi)(\pa_t\Gamma^\mu\phi)
-\sum_{i,j=1}^2\pa_t(\pa_i\phi\pa_j\phi)(\pa_j\Gamma^\mu\phi)(\pa_i\Gamma^\mu\phi)\\
&-\sum_{i=1}^2\pa_i(|\nabla\phi|^2)(\pa_i\Gamma^\mu\phi)(\pa_t\Gamma^\mu\phi)+\frac{1}{2}\pa_t(|\nabla\phi|^2)|\nabla \Gamma^\mu\phi|^2\Big\} dx,\\
L^\mu_4(t)=&(t+2M+2)^{-1}\int_{\Bbb R^2}\tilde\si^{-1/2}e^{\varphi}
\tilde\chi'\Big\{\frac{r}{t+2M+2}\big((\pa_t\Gamma^\mu\phi)^2
+|\nabla\Gamma^\mu\phi|^2\big)-2\sum_{i=1}^2\omega_i\pa_i\phi(\pa_t\Gamma^\mu\phi)^2\\
&+2\sum_{i=1}^2\omega_i(\pa_t\Gamma^\mu\phi)(\pa_i\Gamma^\mu\phi)
+4\sum_{i=1}^2\omega_i\pa_t\phi(\pa_t\Gamma^\mu\phi)
(\pa_i\Gamma^\mu\phi)+2r(t+2M+2)^{-1}\pa_t\phi|\nabla\Gamma^\mu\phi|^2\\
&-2\sum_{i,j=1}^2\omega_i\pa_i\phi\pa_j\phi(\pa_j\Gamma^\mu\phi)(\pa_t\Gamma^\mu\phi)
-r(t+2M+2)^{-1}\pa_i\phi\pa_j\phi(\pa_j\Gamma^\mu\phi)(\pa_i\Gamma^\mu\phi)\\
&+2\sum_{i=1}^2\omega_i|\nabla\phi|^2(\pa_t\Gamma^\mu\phi)(\pa_i\Gamma^\mu\phi)
+r(t+2M+2)^{-1}|\nabla\phi|^2|\nabla\Gamma^\mu\phi|^2\Big\} dx,\\
L^\mu_5(t)=&\frac{1}{2}\int_{\Bbb R^2}(|x|-t-M)\tilde\si^{-5/2}e^{\varphi}
\chi\Big\{\sum_{i=1}^2\omega_i\pa_i \phi(\pa_t\Gamma^\mu\phi)^2-2\sum_{i=1}^2\omega_i\pa_t\phi(\pa_i\Gamma^\mu\phi)(\pa_t\Gamma^\mu\phi)\\
-&\sum_{i=1}^2\pa_t\phi(\pa_i\Gamma^\mu\phi)^2+\sum_{i,j=1}^2\omega_i\pa_i\phi\pa_j\phi(\pa_j\Gamma^\mu\phi)
(\pa_t\Gamma^\mu\phi)+\frac{1}{2}\sum_{i,j=1}^2\pa_i\phi\pa_j\phi(\pa_i\Gamma^\mu\phi)(\pa_j\Gamma^\mu\phi)\\
-&\sum_{i=1}^2\omega_i|\nabla\phi|^2(\pa_i\Gamma^\mu\phi)(\pa_t\Gamma^\mu\phi)
-\frac{1}{2}|\nabla\phi|^2|\nabla\Gamma^\mu\phi|^2\Big\} dx.
\end{align*}

By $\pa_t\dot\phi=-\dot\theta+F_2(\dot\xi,2\xi_a+\dot\xi)$ and $\nabla\dot\phi=\dot w$ for $|x|\geq M+1$, then
by Lemma 2.3 and Lemma 3.5 we have
$$
|\pa\Gamma^\nu\phi(t,x)|\leq C\ve\si_-(t,x)^{-1/2}\si(x)^{-1/2}\quad\text{for}\quad |\nu|\leq \la-2\quad\text{and}\quad |x|\geq M+1.
\eqno{(4.17)}
$$

Similarly,
$$
|\pa^2\Gamma^\nu \phi(t,x)|\leq C\ve\si_-(t,x)^{-1}\si(x)^{-1/2}\quad\text{for}\quad |\nu|\leq \la-3\quad\text{and}\quad |x|\geq M+1.
\eqno{(4.18)}
$$

In addition, it follows Lemma 3.4 that
$$
|\Gamma^\nu\phi(t,x)|\leq C\ve\si(x)^{-1/2}\si_-(t,x)\quad\text{for}\quad |\nu|\leq \la-2\quad\text{and}\quad |x|\geq M+1.
\eqno{(4.19)}
$$

Therefore, we derive from (4.16)-(4.19) that
\begin{align*}
&\widetilde E_k(t)+\sum_{|\mu|\leq k}\int_{1/\ve}^t\int_{\Bbb R^2}\tilde\si^{-1/2}e^\varphi\chi\varphi'|Z\Gamma^\mu\phi|^2 dx d\tau\\
&\qquad -\sum_{|\mu|\leq k}\int_{1/\ve}^t\int_{\Bbb R^2}(|x|-t-M)\tilde\si^{-5/2}e^{\varphi}\chi|Z\Gamma^\mu\phi|^2dxd\tau\\
&\leq C\Big(\widetilde E_k(\frac{1}{\ve})+\sum_{|\mu|\leq k}\sum_{i=1}^5\int_{1/\ve}^t|L^\mu_i(\tau)|d\tau\Big).\tag{4.20}
\end{align*}

Next we deal with each term $L^\mu_i (1\le i\le 5)$ in the right hand side of (4.20).
In this process, we will often use the fact that $\tilde\si(t,x)$ is equivalent to $\si_-(t,x)$.

First, we take $g_i^{0i}=g_i^{i0}=-1$, $g_0^{ii}=2$, $g_{ij}^{ij}=g_{ij}^{ji}=-\ds\frac{1}{2}$, $g_{ii}^{jj}=1$ for $i,j=1,2$, and the others are 0, then we have by Lemma 3.1 and Lemma 3.2 together with (4.17)-(4.19)
\begin{align*}
&|L_1^\mu (t)|\leq \int_{D_+(t)}\tilde\si^{-1/2}e^\varphi\chi|\Gamma^\mu\big(\sum_{i,j,k=0}^2g_i^{jk}\pa_i\phi\pa_{jk}^2\phi\big)
-\sum_{i,j,k=0}^2g_i^{jk}\pa_i\phi\pa_{jk}^2\Gamma^\mu\phi||\pa_t\Gamma^\mu\phi|{d}x\\
&\quad +\sum_{\nu\leq \mu}\int_{D_-(t)}\tilde\si^{-1/2}e^\varphi\chi|C_{\nu}\Gamma^\nu\big(\sum_{i,j,k=0}^2g_i^{jk}\pa_i\phi\pa_{jk}^2\phi\big)
-\sum_{i,j,k=0}^2g_i^{jk}\pa_i\phi\pa_{jk}^2\Gamma^\mu\phi||\pa_t\Gamma^\mu\phi|{d}x\\
&\quad +\sum_{\nu\leq\mu}\int_{\Bbb{R}^2}\tilde\si^{-1/2}e^\varphi\chi|C_{\nu}\Gamma^\nu\big(\sum_{i,j,k,l=0}^2g_{ij}^{kl}\pa_i\phi\pa_j\phi
\pa_{kl}^2\phi\big)-\sum_{i,j,k,l=0}^2g_{ij}^{kl}\pa_i\phi\pa_j\phi\pa_{kl}^2\Gamma^\mu\phi||\pa_t\Gamma^\mu\phi|{d}x\\
&\quad +\sum_{\nu<\mu}\int_{D_+(t)}\tilde\si^{-1/2}e^\varphi\chi|C_{\nu}\Gamma^\nu
\big(\sum_{i,j,k=0}^2g_i^{jk}\pa_i\phi\pa_{jk}^2\phi\big)||\pa_t\Gamma^\mu\phi|dx\\
&\leq C\Big\{\langle t\rangle^{-1}\int_{D_+(t)}\tilde\si^{-1/2}e^\varphi\chi\big(\si_-\sum_{|\nu|\leq [\frac{k+1}{2}]}|\pa \Gamma^\nu \phi|+\sum_{|\nu|\leq [\frac{k+1}{2}]+1}|\Gamma^\nu \phi|\big)\sum_{|\nu|\leq k}|\pa \Gamma^\nu \phi|^2{d}x\\
&\quad +\int_{D_-(t)}\tilde\si^{-1/2}e^\varphi\chi\sum_{|\nu_1|\leq[\frac{k+1}{2}]}|\pa \Gamma^{\nu_1}\phi|\sum_{|\nu_2|\leq k}|\pa \Gamma^{\nu_2} \phi|^2{d}x\\
&\quad +\int_{\Bbb{R}^2}\tilde\si^{-1/2}e^\varphi\chi\sum_{|\nu_1|\leq [\frac{k+1}{2}]}|\pa \Gamma^{\nu_1}\phi|^2\sum_{|\nu_2|\leq k}|\pa \Gamma^{\nu_2}\phi|^2{d}x\\
&\quad +\int_{D_+(t)}\tilde\si^{-1/2}e^\varphi\chi\sum_{|\nu_1|\leq[\frac{k}{2}]}|\pa \Gamma^{\nu_1}\phi|\sum_{|\nu_2|\leq k}|Z\Gamma^{\nu_2}\phi|\sum_{|\nu_3|\leq k}|\pa \Gamma^{\nu_3}\phi|{d}x\Big\}\tag{4.21}\\
&\leq C\ve \sum_{|\nu|\leq k}\int_{\Bbb{R}^2}\tilde\si^{-1/2}e^\varphi\chi{\varphi'}|Z\Gamma^\nu\phi|^2{d}x+C\ve \langle t\rangle^{-1}\sum_{|\nu|\leq k}\int_{\Bbb{R}^2}e^\varphi\chi|\pa\Gamma^\nu\phi|^2{d}x,\tag{4.22}
\end{align*}
here we give some explanations for the derivation process from (4.21) to (4.22), for example,
in order to treat the last term in (4.21), one can make use of (4.17) to get
\begin{align*}
&\int_{D_+(t)}\tilde\si^{-1/2}e^\varphi\chi\sum_{|\nu_1|\leq[\frac{k}{2}]}|\pa \Gamma^{\nu_1}\phi|\sum_{|\nu_2|\leq k}|Z\Gamma^{\nu_2}\phi|\sum_{|\nu_3|\leq k}|\pa \Gamma^{\nu_3}\phi|{d}x\\
&\le C\ve\int_{D_+(t)}{\t\si}^{-\f12}e^{\vp}\chi\si_-^{-\f12}\si^{-\f12}\sum_{|\nu_2|\leq k}|Z\Gamma^{\nu_2}\phi|\sum_{|\nu_3|\leq k}|\pa \Gamma^{\nu_3}\phi|{d}x\\
&\leq C\ve \sum_{|\nu|\leq k}\int_{\Bbb{R}^2}\tilde\si^{-1/2}e^\varphi\chi{\varphi'}|Z\Gamma^\nu\phi|^2{d}x+C\ve \langle t\rangle^{-1}\sum_{|\nu|\leq k}\int_{\Bbb{R}^2}e^\varphi\chi|\pa\Gamma^\nu\phi|^2{d}x\\
&\qquad \quad\qquad \quad\qquad \quad\qquad \quad\text{($\si_-$
is equivalent to $\t\si$ in $D_+(t)$)}.
\end{align*}

Analogously,
\begin{align*}
|L_2^\mu(t)|&\leq C\Big\{\int_{D_+(t)}\tilde\si^{-1/2}e^\varphi\chi{\varphi'}|\sum_{i=1}^2\omega_i Z_i\phi(\pa_t\Gamma^\mu \phi)^2-\pa_t \phi\sum_{i=1}^2(Z_i\Gamma^\mu\phi)^2|{d}x\\
&\quad +\int_{D_-(t)}\tilde\si^{-1/2}e^\varphi\chi{\varphi'}|\pa\phi||\pa\Gamma^\mu \phi|^2{d}x+\int_{\Bbb{R}^2}\tilde\si^{-1/2}e^\varphi\chi{\varphi'}|\pa\phi|^2|\pa\Gamma^\mu\phi|^2{d}x\Big\}\\
&\leq C\ve\langle t\rangle^{-1/2}\int_{\Bbb{R}^2}\tilde\si^{-1/2}e^\varphi\chi{\varphi'}|Z\Gamma^\mu\phi|^2{d}x+C\ve\langle t\rangle^{-1}\int_{\Bbb{R}^2}e^\varphi\chi|\pa\Gamma^\mu\phi|^2{d}x,\tag{4.23}\\
|L_3^\mu(t)|&\leq C\Big\{\int_{D_+(t)}\tilde\si^{-1/2}e^\varphi\chi\big(|\pa\Gamma^\mu\phi|^2|Z\pa\phi|
+|(\pa\Gamma^\mu\phi)(Z\Gamma^\mu\phi)(\pa^2\phi)|\big) dx\\
&\quad +\int_{D_-(t)}\tilde\si^{-1/2}e^\varphi\chi|(\pa^2\phi)(\pa\Gamma^\mu\phi)^2| dx
+\int_{\Bbb R^2}\tilde\si^{-1/2}e^\varphi\chi|\pa\phi(\pa^2\phi)(\pa\Gamma^\mu\phi)^2| dx\Big\}\\
&\leq C\ve\int_{\Bbb{R}^2}\tilde\si^{-1/2}e^\varphi\chi{\varphi}'|Z\Gamma^\mu\phi|^2{d}x+C\ve\langle t\rangle^{-1}\int_{\Bbb{R}^2}e^\varphi\chi|\pa\Gamma^\mu\phi|^2{d}x,\tag{4.24}\\
|L_4^\mu(t)|&\leq C\langle t\rangle^{-1}\int_{\{|x|\geq M+1\}}e^\varphi|\pa\Gamma^\mu\phi|^2{d}x,\tag{4.25}\\
|L_5^\mu(t)|&\leq C\Big\{-\int_{D_+(t)}(|x|-t-M)\tilde\si^{-5/2}e^\varphi\chi|\sum_{i=1}^2\omega_i Z_i\phi(\pa_t\Gamma^\mu\phi)^2-\pa_t\phi\sum_{i=1}^2(Z_i\Gamma^\mu\phi)^2|{d}x\\
&\quad +\int_{D_-(t)}\tilde\si^{-3/2}e^\varphi\chi|\pa\phi||\pa\Gamma^\mu\phi|^2{d}x+\int_{\Bbb{R}^2}\tilde\si^{-3/2}e^\varphi\chi|\pa \phi|^2|\pa\Gamma^\mu\phi|^2{d}x\Big\}\\
&\leq -C\ve\langle t\rangle^{-1/2}\int_{\Bbb{R}^2}(|x|-t-M)\tilde\si^{-5/2}e^\varphi\chi|Z\Gamma^\mu\phi|^2{d}x+C\ve\langle t\rangle^{-1}\int_{{R}^2}e^\varphi\chi|\pa\Gamma^\mu\phi|^2{d}x.\tag{4.26}
\end{align*}

Substituting (4.22)-(4.26) into (4.20) yields
\begin{align*}
&\widetilde E_k(t)+\sum_{|\mu|\leq k}\int_{1/\ve}^t\int_{\Bbb R^2}
\tilde\si^{-1/2}e^\varphi\chi\varphi'|Z\Gamma^\mu\phi|^2dx d\tau\\
&\qquad -\sum_{|\mu|\leq k}\int_{1/\ve}^t\int_{\Bbb R^2}(|x|-t-M)\tilde\si^{-5/2}e^{\varphi}\chi|Z\Gamma^\mu\phi|^2 dx d\tau\\
&\leq C\widetilde E_k(\frac{1}{\ve})+C\int_{1/\ve}^t\langle \tau\rangle^{-1}\sum_{|\mu|\leq k}\int_{\{|x|\geq M+1\}}|\pa\Gamma^\mu \phi|^2{d}x.\tag{4.27}
\end{align*}

This, together with Lemma 4.1, yields (4.13).\qquad\qquad\qquad\qquad $\square$ \bigskip

From (4.13) and (2.2), we have for $\dot\phi=\phi-\phi_a$
\begin{align*}
&\sum_{|\mu|\leq k}\int_{1/\ve}^t\int_{D_+(\tau)}\si_-(t,x)^{-2}|Z\Gamma^\mu\dot\phi(\tau,x)|^2 dx d\tau\\
&\leq C\Big(\ve^2\ln t+\sum_{|\mu|\leq k}\int_{1/\ve}^t\langle\tau\rangle^{-1}\int_{\{|x|\geq M+1\}}|\Gamma^\mu\pa\dot\phi(\tau,x)|^2dx d\tau\Big).\tag{4.28}
\end{align*}

{\bf Lemma 4.3.} {\it
Assume that $k$ and $\la$ are integers with $[\frac{k+7}{2}]\leq\la\leq k$ and $k\geq7$,
$(\dot\theta, \dot w, \dot z)$ is a smooth solution of
(2.4) for $(t,x)\in [1/\ve, T]\times \Bbb R^2$.
If $E_\la(t)\leq\ve^2$, then one can find a positive number $C$ independent of $\ve$ and $T$ such that for
$t\in[1/\ve,T]$ and sufficient small $\ve>0$,
$E_k(t)\leq C\ve^2\langle t\rangle^{C\ve}$ holds.
}

{\bf Proof.}
First, we come to estimate $H_j^\mu$ in (4.5) when $|\mu|\leq k$ and $j=1,\cdots,13$. It is easy to conclude that
$$
\sum_{|\mu|\leq k}\sum_{i=1}^{13}|H_i^\mu(t)|_-\leq C\ve\langle t\rangle^{-1}E_k(t)
+C\ve^2\langle t\rangle^{-5/2}E_k(t)^{1/2},\eqno{(4.29)}
$$
here the estimate on $H_{13}^\mu$ can be obtained as in (4.11).

We now continue to use the analogous notations as in (3.21), namely,
$$
\hat{h}_0^\mu=\sum_{0<\nu\leq\mu}C_{\mu\nu}\big(I_1^{\mu\nu}+I_3^{\mu\nu}\big)
+\sum_{0\leq \nu\leq\mu}C_{\mu\nu}I_2^{\mu\nu}+K^\mu,
$$
where
\begin{align*}
K^\mu=&\sum_{0<\nu\leq\mu}C_{\mu\nu}\big(-\Gamma^\nu w\cdot\nabla\Gamma^{\mu-\nu}F_2(\dot\xi,2\xi_a+\dot\xi)+\Gamma^\nu F_1(\xi)\nabla\cdot\Gamma^{\mu-\nu}\dot w\big)\\
&+\sum_{0\leq \nu\leq\mu}C_{\mu\nu}\big(-\Gamma^\nu \dot w\cdot\nabla\Gamma^{\mu-\nu}F_1(\xi_a)+\Gamma^\nu F_2(\dot\xi,2\xi_a+\dot\xi)\nabla\cdot\Gamma^{\mu-\nu}w_a\big).
\end{align*}

It follows from the similar analysis of (3.22) and (3.24) together with (4.28) that
\begin{align*}
&\sum_{0<\nu\leq\mu}C_{\mu\nu}\|I_1^{\mu\nu}(t)\|_+\leq C\ve \langle t\rangle^{-3/2}E^{1/2}_k(t),\qquad\qquad\qquad\qquad\qquad\qquad\qquad\tag{4.30}\\
&\ve^{-1}\sum_{0\leq \nu\leq \mu}C_{\mu\nu}\int_{1/\ve}^t\langle\tau\rangle\|I_2^{\mu\nu}(\tau)\|_+^2d\tau\\
&\quad \leq C\ve^{-1}\int_{1/\ve}^t\langle\tau\rangle\int_{D_+(\tau)}\sum_{|\nu_1|+|\nu_2|\leq k}\Big(|Z\Gamma^{\nu_1} \dot\phi|^2|\pa^2\Gamma^{\nu_2} \phi_a|^2+|\pa\Gamma^{\nu_1} \dot\phi|^2|Z\pa\Gamma^{\nu_2} \phi_a|^2\Big) dx d\tau\\
&\quad \leq C\ve\sum_{|\nu|\leq k}\int_{1/\ve}^t\int_{D_+(\tau)}\si_-(t,x)^{-2}|Z\Gamma^\nu\dot\phi(t,x)|^2dx d\tau+C\ve\int_{1/\ve}^t\langle\tau\rangle^{-1}E_k(\tau) d\tau\\
&\quad \leq C\ve^3\ln t+C\ve\int_{1/\ve}^t\langle\tau\rangle^{-1}E_k(\tau) d\tau,\tag{4.31}\\
&\ve^{-1}\sum_{0< \nu\leq\mu}C_{\mu\nu}\int_{1/\ve}^t\langle\tau\rangle\|I_3^{\mu\nu}(\tau)\|_+^2d\tau\\
&\quad \leq C\ve^{-1}\int_{1/\ve}^t\langle\tau\rangle\int_{D_+(\tau)}\sum_{|\nu_1|+|\nu_2|\leq k,|\nu_2|\leq k-1}\Big(|Z\Gamma^{\nu_1}\dot\phi|^2|\pa^2\Gamma^{\nu_2} \dot\phi|^2+|\pa\Gamma^{\nu_1}\dot\phi|^2|Z\pa\Gamma^{\nu_2}\dot\phi|^2\Big) dx d\tau\\
&\quad \leq C\ve\sum_{|\nu|\leq k-1} \int_{1/\ve}^t\langle\tau\rangle^{-2}\int_{D_+(\tau)}|\si_-\pa\Gamma^\nu\pa\dot\phi|^2{d}x d\tau\\
&\qquad +C\ve\sum_{|\nu|\leq k}\int_{1/\ve}^t\int_{D_+(\tau)}\si_-(x,t)^{-2}|Z\Gamma^\nu\dot\phi(t,x)|^2dxd\tau
+C\ve\int_{1/\ve}^t\langle\tau\rangle^{-1}E_k(\tau)d\tau\\
&\quad \leq C\ve^3\ln t+C\ve\int_{1/\ve}^t\langle\tau\rangle^{-1}E_k(\tau)d\tau.\tag{4.32}
\end{align*}

By using Lemma 2.2 and Lemma 2.3, we can get
$$
|K^\mu(t)|_+\leq C\ve^2\langle t\rangle^{-1}\sum_{|\nu|\leq k}\big(|\Gamma^\nu \dot\theta(t)|_+
+|\Gamma^\nu\dot w(t)|_+\big).\eqno{(4.33)}
$$

Hence, we have
$$
\int_{1/\ve}^t\sum_{|\mu|\leq k}\sum_{i=1}^{6}|H_i^\mu(\tau)|_+d\tau\leq C\ve^3\ln t+C\ve\int_{1/\ve}^t\langle\tau\rangle^{-1}E_k(\tau) d\tau,
$$
and similarly,
$$
\int_{1/\ve}^t\sum_{|\mu|\leq k}\sum_{i=7}^{12}|H_i^\mu(t)|_+ d\tau\leq C\ve^3\ln t+C\ve\int_{1/\ve}^t\langle\tau\rangle^{-1}E_k(\tau) d\tau.
$$

Second, we deal with the terms $\hat H^\mu_{aj}$ and $\hat H^\mu_{j}$ in (4.5) for $j=1,2$ and $|\mu|\leq k$.

Due to $\text{supp}\dot z\subset\{|x|\leq M+1\}$,
it is obvious that for $j=1,2$
$$|\hat H_{aj}^\mu|\leq C\ve\langle t\rangle^{-5/2}E_k^{1/2}(t)\big(E_k^{1/2}(t)+\ve\big).\eqno{(4.34)}$$

To treat the terms $\hat H_{j}^\mu$ ($j=1,2$), we will use the analogous method in $\S 5$ of [9] (see pages 101 of [9]).
For this end, we set $\la_a^{\mu\nu}=\Gamma^\nu \dot w\cdot\nabla\Gamma^{\mu-\nu}z(0,x)$ and $\la^{\mu\nu}=\Gamma^\nu \dot w\cdot\nabla\Gamma^{\mu-\nu}\dot z$.

If $\Gamma^\nu=\pa\Gamma^d$ with $|d|=|\nu|-1$, then
$$
\|\la_a^{\mu\nu}(t)\|\leq C\ve\langle t\rangle^{-1}\|(\si_-\pa\Gamma^d\dot w)(t)\|\leq C\ve\langle t\rangle^{-1}\big(E_k^{1/2}(t)+\ve^2\langle t\rangle^{-3/2}\big)\eqno{(4.35)}
$$
and
$$
\|\la^{\mu\nu}(t)\|\leq\big(\sup_{|x|\leq M+1}|\pa\Gamma^d\dot w(t,x)|\big)\|\nabla\Gamma^{\mu-\nu}\dot z(t)\|\leq C\ve\langle t\rangle^{-1}E_k^{1/2}(t)\quad\text{for}\quad |\nu|\leq|\mu-\nu|,\eqno{(4.36)}
$$
\begin{align*}
&\|\la^{\mu\nu}(t)\|\leq C\big(\sup\nabla\Gamma^{\mu-\nu} \dot z(t,x)\big)\langle t\rangle^{-1}\|\si_-\pa\Gamma^d\dot w(t)\|\leq C\ve\langle t\rangle^{-1}\big(E_k^{1/2}(t)+\ve^2\langle t\rangle^{-3/2}\big)\\
&\qquad \quad\qquad \quad\qquad \quad\text{for}\quad |\nu|>|\mu-\nu|.\tag{4.37}
\end{align*}

If $\Gamma^\nu=X^{|\nu|}$ with $|\nu|\leq k-1$, due to $\Gamma^\nu \dot w(t,x)=\big(\Gamma^\nu W(t,r)\big)\ds\frac{x}{r}$
holds for some smooth function $W(t,r)$, then we have
$$
\sup_{|x|\leq M+1}|\Gamma^\nu\dot w(t,x)|\leq C\langle t\rangle^{-1}\|\si_-\nabla\cdot\Gamma^\nu\dot w(t,x)\|\leq C\langle t\rangle^{-1}Q_{|\nu|+1}(t).\eqno{(4.38)}
$$

Hence, we derive from (4.38) that for $|\nu|\leq k-1$
$$
\|\la_a^{\mu\nu}(t)\|\leq C\ve\langle t\rangle^{-1}Q_k(t)\leq C\ve\langle t\rangle^{-1}\big(E_k^{1/2}(t)+\ve^2\langle t\rangle^{-3/2}\big),\eqno{(4.39)}
$$
$$
\|\la^{\mu\nu}(t)\|\leq C\langle t\rangle^{-1}Q_{|\nu|+1}(t)E_k^{1/2}(t)\leq C\ve\langle t\rangle^{-1}E_k^{1/2}(t)\quad\text{for}\quad |\nu|\leq|\mu-\nu|,\eqno{(4.40)}
$$
and
$$
\|\la^{\mu\nu}(t)\|\leq C\ve\langle t\rangle^{-1}\big(E_k^{1/2}(t)
+\ve^2\langle t\rangle^{-3/2}\big)\quad\text{for}\quad |\nu|>|\mu-\nu|.\eqno{(4.41)}
$$

If $\Gamma^\nu=X^k$, then exactly as in the proof of (5.10) in [9]
(here we will apply Lemma 2.5) together with  the related estimates (4.35)-(4.41),
we can obtain
$$
|\hat H_{j}^\mu(t)-\frac{d}{dt}G_j(t)|\leq C\ve\langle t\rangle^{-1}\big(E_k(t)
+\ve^2\langle t\rangle^{-3/2}E^{1/2}_k(t)+\ve^5\langle t\rangle^{-3}\big)\qquad \text{for $j=1,2$},\eqno{(4.42)}
$$
where $|G_j(t)|\leq C\big(\ve E_k(t)+\ve^3\langle t\rangle^{-3/2}E^{1/2}_k(t)+\ve^6\langle t\rangle^{-3}\big)$.

Combining (4.34) with (4.42) yields
$$
\sum_{|\mu|\leq k}\sum_{j=1,2}(\hat H^\mu_{aj}(t)+\hat H^\mu_{j}(t))=\frac{ d}{ dt}\t G_1(t)+\t G_2(t),\eqno{(4.43)}
$$
where $|\t G_1(t)|\leq C\big(\ve E_k(t)+\ve^3\langle t\rangle^{-3/2}E^{1/2}_k(t)+\ve^6\langle t\rangle^{-3}\big)$, and $|\t G_2(t)|\leq C\ve\langle t\rangle^{-1}\big(E_k(t)+\ve^5\langle t\rangle^{-3}$ $+\ve\langle t\rangle^{-3/2}E^{1/2}_k(t)\big)$.

Third, according to the expressions of $A_j$ for $j=0, 1, 2$, it is easily known that
$$
\sum_{j=0}^2|\pa_j A_j(t)|_-\leq C\big(|\nabla\theta(t)|_-+|\nabla\cdot w(t)|_-+|w\cdot\nabla z(t)|_-\big)\leq C\ve\langle t\rangle^{-3/2},
$$
and then
$$|\ds\sum_{j=0}^2\langle(\pa_j A_j)\zeta^\mu,\zeta^\mu\rangle_-(t)|\leq C\ve\langle t\rangle^{-3/2} E_k(t).\eqno{(4.44)}$$

For the case of $|x|\geq \ds\frac{t}{2}
+M+1$, we know $z(t,x)=0$ by Lemma 2.6 and
\begin{align*}
&\sum_{j=0}^2\langle(\pa_j A_j)\zeta^\mu,\zeta^\mu\rangle_+(t)\\
=&\langle(\nabla\cdot w)\Gamma^\mu\dot\theta-\nabla\theta\cdot\Gamma^\mu\dot w,\Gamma^\mu
\dot\theta\rangle_+(t)+\langle(\nabla\cdot w)\Gamma^\mu \dot w-(\Gamma^\mu\dot\theta)\nabla\theta,\Gamma^\mu\dot w\rangle_+(t)\\
=&\langle \pa_t\nabla\phi\cdot\Gamma^\mu\nabla\dot\phi-\triangle\phi\Gamma^\mu\pa_t\dot\phi,\Gamma^\mu
\dot\theta\rangle_+(t)+\langle(\triangle\phi)\Gamma^\mu\nabla\dot\phi-(\Gamma^\mu\pa_t\dot\phi)\nabla\pa_t\phi,\Gamma^\mu \dot w\rangle_+(t)\\
&+\langle (\nabla\cdot w)\Gamma^\mu F_2(\dot\xi,2\xi_a+\dot\xi)
-\nabla F_1(\xi)\cdot\Gamma^\mu \dot w,\Gamma^\mu\dot\theta\rangle_+(t)\\
&+\langle(\Gamma^\mu F_2(\dot\xi,2\xi_a
+\dot\xi)-\Gamma^\mu\dot\theta)\nabla F_1(\xi)-\Gamma^\mu F_2(\dot\xi,2\xi_a+\dot\xi)\nabla\theta,
\Gamma^\mu\dot w\rangle_+(t).\tag{4.45}
\end{align*}

Notice that
\begin{align*}
&\int_{1/\ve}^t\langle \pa_t\nabla\phi\cdot\Gamma^\mu\nabla\dot\phi-\triangle\phi\Gamma^\mu\pa_t\dot\phi,\Gamma^\mu
\dot\theta\rangle_+(\tau)d\tau\\
&\quad \leq C\sum_{l=1}^2\sum_{\nu\leq\mu}\int_{1/\ve}^t|\langle\pa_t\pa_l\phi\pa_l\Gamma^\nu \dot\phi-\pa_l^2\phi\pa_t\Gamma^\nu\dot\phi,\Gamma^\mu\dot\theta\rangle_+(\tau)|
d\tau\\
&\quad \leq C\ve\sum_{|\nu|\leq k}\int_{1/\ve}^t\int_{D_+(\tau)}\si_-^{-2}|Z\Gamma^\nu\dot\phi|^2dxd\tau+C\ve\int_{1/\ve}^t
\langle\tau\rangle^{-1}E_k(\tau)d\tau\\
&\quad \leq C\ve^3\ln t+C\ve\int_{1/\ve}^t\langle\tau\rangle^{-1}E_k(\tau)d\tau,\tag{4.46}
\end{align*}
and in the same way, one has
$$
\int_{1/\ve}^t\langle(\triangle\phi)\Gamma^\mu\nabla\dot\phi-(\Gamma^\mu\pa_t\dot\phi)\nabla\pa_t\phi,
\Gamma^\mu\dot w\rangle_+(\tau)d\tau
\leq C\ve^3\ln t+C\ve\int_{1/\ve}^t\langle\tau\rangle^{-1}E_k(\tau)d\tau\eqno{(4.47)}
$$
and
\begin{align*}
&\langle (\nabla\cdot w)\Gamma^\mu F_2(\dot\xi,2\xi_a+\dot\xi)-\nabla F_1(\xi)\cdot\Gamma^\mu\dot w,
\Gamma^\mu\dot\theta\rangle_+(t)\\
&\quad+\langle(\Gamma^\mu F_2(\dot\xi,2\xi_a+\dot\xi)-\Gamma^\mu\dot\theta)\nabla F_1(\xi)-\Gamma^\mu F_2(\dot\xi,2\xi_a+\dot\xi)\nabla\theta,\Gamma^\mu \dot w\rangle_+(t)\\
&\leq C\ve^2\langle t\rangle^{-1}E_k(t).\tag{4.48}
\end{align*}

Substituting (4.46)-(4.48) into (4.45) and further combining with (4.44) yield
$$
\int_{1/\ve}^t\sum_{j=0}^2\langle(\pa_j A_j)\zeta^\mu,\zeta^\mu\rangle(\tau) d\tau
\leq C\ve^3\ln t+C\ve\int_{1/\ve}^t\langle\tau\rangle^{-1}E_k(\tau)d\tau.\eqno{(4.49)}
$$

Based on the estimates of the above three steps, if we set $N(t)=\ds\sum_{|\mu|\leq k}\langle A_0\zeta^\mu,
\zeta^\mu\rangle(t)-2\t G_1(t)$, then
$$
N(t)\leq C\ve^3\ln t+C\ve\int_{1/\ve}^t\langle\tau\rangle^{-1}E_k(\tau)d\tau
+C\ve^2\int_{1/\ve}^t\langle\tau\rangle^{-5/2}E_k^{1/2}(\tau)d\tau.\eqno{(4.50)}
$$

Note that $\tilde N(t)=N(t)+\ve^4\langle t\rangle^{-3}$ is equivalent to $E_k(t)+\ve^4\langle t\rangle^{-3}$,
then one can derive from (4.50)
$$
\tilde N(t)\leq C\ve^3\ln t+C\ve\int_{1/\ve}^t\langle\tau\rangle^{-1}\tilde N(\tau)d\tau
+C\ve^2\int_{1/\ve}^t\langle\tau\rangle^{-5/2}\tilde N^{1/2}(\tau) d\tau.\eqno{(4.51)}
$$

Set $g(t)=C\ve^3\ln t+C\ve\int_{1/\ve}^t\langle\tau\rangle^{-1}\tilde N(\tau) d\tau+C\ve^2\int_{1/\ve}^t\langle\tau\rangle^{-5/2}\tilde N^{1/2}(\tau) d\tau$, it follows from (4.51) that
$$
g'(t)\leq C\ve\langle t\rangle^{-1}\big(g(t)+\ve^2\big)+C\ve^2\langle t\rangle^{-5/2}g(t)^{1/2},
$$
furthermore, if we set $\tilde g(t)=g(t)+\ve^2$, then
$$
 \tilde g'(t)\leq C\ve\langle t\rangle^{-1}\tilde g(t)+C\ve^2\langle t\rangle^{-5/2}\tilde g(t)^{1/2},
$$
thus, $g(t)\leq \tilde g(t)\leq C\ve^{2}\langle t\rangle^{C\ve}$ and further $E_k(t)\leq C\ve^2\langle t\rangle^{C\ve}$.
\qquad $\square$\bigskip

Based on Lemma 4.3, we next derive the uniform energy estimate on the solution $(\dot\th, \dot w, \dot z)$ of (2.4).

{\bf Lemma 4.4.} {\it
For a fixed integer $\la$ with $\la\geq 9$, it is assumed that $(\dot\theta, \dot w, \dot z)$ is a smooth solution of
(2.4) for $(t,x)\in [1/\ve, T]\times\Bbb R^2$.
If $E_\la(t)\leq\ve^2$ for $t\in[1/\ve,T]$, then $E_\la(t)\leq \frac{1}{2}\ve^2$ holds for small $\ve>0$.
}
\bigskip

{\bf Proof.}
Similar to the proof of Lemma 4.3, we will divide the whole proof procedure into three steps.
In the following, we always assume $|\mu|\leq\la$ and apply the same notations in (4.5).

First, by Lemma 2.2 and the assumption of $E_\la(t)\leq\ve^2$ for $t\in[1/\ve,T]$, we can get
$$
|H_j^\mu(t)|_-\leq C\ve^3\langle t\rangle^{-5/2},\quad j\in\{1,2,4,5,7,8,10,11\}.\eqno{(4.52)}
$$

Since $|\Gamma^\nu\dot w(t)|_-+|\Gamma^\nu\dot\theta(t)|_-\leq C\langle t\rangle^{-1/2}\widetilde Q_{\mu+2}(t)\leq C\ve\langle t\rangle^{-1/2+\delta}$ holds by Lemma 4.3 for sufficient small positive number $\delta$, we have
$$
|H_j^\mu(t)|_-\leq C\ve^3\langle t\rangle^{-3/2+\delta},\quad j\in\{3,6,9,12\}.\eqno{(4.53)}
$$

In addition, from (4.9) and (4.10) we can obtain $\|\hat f^\mu_{13}\|\leq C\ve^2\langle t\rangle^{-1}$ and further
$$
|H_{13}^\mu(t)|\leq C\ve^2\langle t\rangle^{-1}\sup_{|x|\leq M+1}|\Gamma^b\dot w|\leq C\ve^3\langle t\rangle^{-3/2+\delta}.\eqno{(4.54)}
$$

On the other hand, it follows (3.21) and (3.23)-(3.27) that
$$
\sum_{j=1}^{12} |H_j^\mu(t)|_+\leq C\ve^3\langle t\rangle^{-3/2+\delta}.\eqno{(4.55)}
$$

Second, it is easy to get
$$
|\hat H^\mu_{aj}(t)|\leq C\ve^3\langle t\rangle^{-5/2}\quad\text{for}\quad j=1,2.\eqno{(4.56)}
$$

Lemma 2.5 gives $|\Gamma^\nu\dot w(t,x)|\leq C\ve\langle t\rangle^{-3/2+\delta}$
for $|x|\leq M+1$ and $|\nu|\leq\mu$, and hence
$$
|\hat H^\mu_{j}(t)|\leq C\ve^3\langle t\rangle^{-3/2+\delta}\quad\text{for}\quad j=1,2.\eqno{(4.57)}
$$

Third, as in the proof of Lemma 4.3, one has
$$
|\sum_{j=0}^2\langle(\pa_j A_j)\zeta^\mu,\zeta^\mu\rangle_-(t)|\leq C\ve^3\langle t\rangle^{-3/2}.\eqno{(4.58)}
$$

Similar to the estimate (3.24) on $I_3^{\mu\nu}$, we have
$$
\|(\pa_t\nabla\phi\cdot\Gamma^\mu\nabla\dot\phi-\triangle\phi\Gamma^\mu\pa_t\dot\phi)(t)\|_+\leq C\ve\langle t\rangle^{-3/2}E_{\la+1}^{1/2}\leq C\ve^2\langle t\rangle^{-3/2+\delta}\eqno{(4.59)}
$$
and
$$
\|(\triangle\phi\Gamma^\mu\nabla\dot\phi-\Gamma^\mu\pa_t\dot\phi\nabla\pa_t\phi)(t)\|_+\leq C\ve\langle t\rangle^{-3/2}E_{\la+1}^{1/2}\leq C\ve^2\langle t\rangle^{-3/2+\delta}.\eqno{(4.60)}
$$

In addition, we can decompose
\begin{align*}
&\langle (\nabla\cdot w)\Gamma^\mu F_2(\dot\xi,2\xi_a+\dot\xi)-\nabla F_1(\xi)\cdot\Gamma^\mu\dot w,\Gamma^\mu\dot\theta\rangle_+(t)\\
&\quad +\langle(\Gamma^\mu F_2(\dot\xi,2\xi_a+\dot\xi)-\Gamma^\mu\dot\theta)\nabla F_1(\xi)-\Gamma^\mu F_2(\dot\xi,2\xi_a+\dot\xi)\nabla\theta,\Gamma^\mu\dot w\rangle_+(t)\\
&\equiv A^\mu_1(t)+A^\mu_2(t)+A^\mu_3(t),\tag{4.61}
\end{align*}
where
\begin{align*}
|A^\mu_1(t)|&=|\langle\triangle\phi\Gamma^\mu F_2(\dot\xi,2\xi_a+\dot\xi)-2\nabla F_1(\xi)\cdot
\Gamma^\mu\nabla \dot\phi,\Gamma^\mu F_2(\dot\xi,2\xi_a+\dot\xi)\rangle_+|\\
&\leq C\ve^3\langle t\rangle^{-3/2}E_\la(t)\leq C\ve^5\langle t\rangle^{-3/2},\tag{4.62}\\
|A^\mu_2(t)|&=|\langle\nabla\pa_t\phi\cdot\Gamma^\mu\nabla\dot\phi-\triangle\phi\Gamma^\mu\pa_t
\dot\phi,\Gamma^\mu F_2(\dot\xi,2\xi_a+\dot\xi)\rangle_+|\leq C\ve^4\langle t\rangle^{-2+\delta},\tag{4.63}\\
|A^\mu_3(t)|&=|2\langle\nabla F_1(\xi)\cdot\Gamma^\mu\nabla\dot\phi,\Gamma^\mu\pa_t\dot\phi\rangle_+|\leq C\ve^2\langle t\rangle^{-3/2}E_\mu(t)\leq C\ve^4\langle t\rangle^{-3/2},\tag{4.64}
\end{align*}
here we point out that the estimate of $\nabla F_1(\xi)$ in $A^\mu_3(t)$ follows from
$$
|\pa_l F_1(\xi)(t)|_+=|\big(\pa_t\phi\pa_l\pa_t\phi-\nabla \phi\cdot\nabla\pa_l \phi\big)(t)+\pa_l\big(F_1(\xi)(\theta-\frac{1}{2}F_1(\xi))\big)(t)|_+
\leq C\ve^2\langle t\rangle^{-3/2}.
$$

Therefore, by substituting (4.62)-(4.64) into (4.61) and further combining with (4.59)-(4.60),
it follows from (4.45) that
$$
|\sum_{j=0}^2\langle(\pa_j A_j)\zeta^\mu,\zeta^\mu\rangle_+(t)|
\leq C\ve^3\langle t\rangle^{-3/2+\delta}.\eqno{(4.65)}
$$

Finally, inserting (4.52)-(4.58) and (4.65) into (4.4) and combining with basic
energy inequalities
(similar to (2.12) and (3.22)) yield
$$\ds\frac{ d}{ dt}E_\mu(t)\leq C\ve^3\langle t\rangle^{-3/2+\delta}.\eqno{(4.66)}$$

In addition, we have $E_\mu(\ds\frac{1}{\ve})\leq C\ve^4$ by Lemma 4.1. This,
together with (4.66),
derives $E_\mu(t)\leq C\ve^3$ when we choose $\delta<\ds\frac{1}{2}$. Therefore,
$E_\mu(t)\leq \ds\f12\ve^2$ holds
for small $\ve>0$.\qquad $\square$

Finally, we start to show Theorem 1.1.
\bigskip

{\bf Proof of Theorem 1.1.} By Lemma 4.1 and Lemma 4.4, we know that (2.4) admits a
global smooth solution $(\dot\th, \dot w,
\dot z)$ in terms of the continuity induction method. Thus, (2.3) has a global smooth
solution $(\th, w, z)$
since the smooth solution of (2.1) exists globally. Therefore, Theorem 1.1 is proved.
\qquad \qquad $\square$

\vskip 0.4 true cm

\end{document}